\documentclass[aos]{imsart} 

\usepackage[reqno]{amsmath}
\RequirePackage{amsfonts,amssymb,amsthm, mathtools}
\RequirePackage[authoryear]{natbib}
\RequirePackage[colorlinks,citecolor=blue,urlcolor=blue]{hyperref}
\RequirePackage{graphicx}

\makeatletter
\usepackage{xcolor}
\usepackage[final]{pdfpages}
\usepackage{ragged2e}
\usepackage{array}
\usepackage{tabularx}
\usepackage{multirow}
\usepackage{rotating}
\usepackage{subcaption}
\usepackage{eurosym}
\usepackage{dsfont}
\usepackage{lmodern}
\usepackage{ulem}
\usepackage{multicol}
\usepackage{hyperref}
\usepackage{etoolbox}
\usepackage{textcomp}
\usepackage{xcolor}
\usepackage{tikz}
\usepackage{tikz-cd}
\usetikzlibrary {positioning}
\usepackage{shadethm}
\usepackage{pdfpages}
\usepackage[utf8]{inputenc} 
\usepackage{csquotes}
\usepackage[T1]{fontenc}
\usepackage[german, english]{babel}
\usepackage{latexsym}
\usepackage{nccmath}
\usepackage{color}
\usepackage{comment}
\usepackage{mathrsfs}
\usepackage{clipboard}
\usepackage{autobreak}
\usepackage{hyperref}
\usepackage{bm}
\usepackage{stackengine,scalerel}
\usepackage{cleveref}
\usepackage{mathtools}
\makeatother

\renewcommand{\thetable}{\arabic{section}.\arabic{table}}
\numberwithin{equation}{section}

\theoremstyle{plain} 
\newtheorem{theorem}{Theorem}[section]
\newtheorem{lemma}[theorem]{Lemma}

\newtheorem{proposition}[theorem]{Proposition}
\newtheorem{corollary}[theorem]{Corollary}
\newtheorem{assumptionletter}{Assumption}
\newshadetheorem{assumption}{Assumption}
\newshadetheorem{property}{Property}
\newshadetheorem{assumptionwn}{Assumption A}

\newshadetheorem{problem}{Optimisation problem}

\theoremstyle{remark}
\newtheorem{definition}[theorem]{Definition}

\newtheorem{example}[theorem]{Example}

\newtheorem{remark}[theorem]{Remark}
\newtheorem{interpretation}[theorem]{Interpretation}

\def\tablename{\footnotesize Table}

\newcommand{\C}{\mathbb{C}} 
\newcommand{\R}{\mathbb{R}} 
\newcommand{\N}{\mathbb{N}} 
\newcommand{\BE}{\mathbb{E}}

\newcommand{\bR}{{R}}

\newcommand{\CX}{{X}}

\newcommand{\CY}{{Y}}

\newcommand{\BA}{\text{\textbf{\textsl{A}}}}
\newcommand{\BB}{\text{\textbf{\textsl{B}}}}
\newcommand{\BFC}{\text{\textbf{\textsl{C}}}}

\newcommand{\obar}[1]{\mkern 1.5mu\overline{\mkern-1.5mu#1\mkern-1.5mu}\mkern 1.5mu}

\newcommand{\oobar}[1]{\obar{\obar{#1}}}
\newcommand{\OBB}{\oobar{\mathbf{B}}}

\newcommand{\OBC}{\text{\textbf{\textsl{C}}}}               
\newcommand{\UBC}{\boldsymbol{\mathsf{C}}}           
\newcommand{\OBE}{\text{\textbf{\textsl{E}}}}               
\newcommand{\UBE}{\boldsymbol{\mathsf{E}}}           
\newcommand{\OBM}{\text{\textbf{\textsl{M}}}}               
\newcommand{\BM}{\boldsymbol{\mathsf{M}}}                
\newcommand{\OBMnu}{\text{\textbf{\textsl{M}}}_m} 
\newcommand{\OBMo}{\text{\textbf{\textsl{M}}}_0} 
\newcommand{\BMnu}{\boldsymbol{\mathsf{M}}_m}    
\newcommand{\BMo}{\boldsymbol{\mathsf{M}}_0} 
\newcommand{\Oeps}{\varepsilon}                  
\newcommand{\Ueps}{\epsilon}                      

\newcommand{\BLAM}{\mathbf{\Lambda}}
\newcommand{\BT}{\mathbf{\Theta}}
\newcommand{\BFE}{\mathbf{E}}
\newcommand{\BS}{\Sigma}

\newcommand{\Astar}{\boldsymbol{\mathcal{A}}}        
\newcommand{\Bstar}{\boldsymbol{\mathcal{B}}}                     
\newcommand{\Cstar}{\boldsymbol{\mathcal{C}}}                 
\newcommand{\Pstar}{\mathcal{P}}        
\newcommand{\Qstar}{\mathcal{Q}}        

\newcommand{\pstar}{\mathfrak{p}}        

\newcommand{\qstar}{\mathfrak{q}}        

\newcommand{\beao}{\begin{eqnarray*}}
\newcommand{\eeao}{\end{eqnarray*}\noindent}
\newcommand{\limm}{\underset{n \rightarrow \infty}{\text{l.i.m.}}}
\newcommand{\limh}{\underset{h \rightarrow 0}{\text{l.i.m.\:}}}

\newcommand{\limhh}{\underset{h \rightarrow 0}{\text{lim\:}}}
\newcommand{\rank}{\text{rk}}

\newcommand{\uint}{\int_{-\infty}^{\infty}}

\newcommand{\rang}{\text{rank}}
\newcommand{\inst}{\:\raisebox{2pt}{\tikz{\draw[-,densely dashed,line width = 0.5 pt](0,0) -- (5mm,0);}}\:}
\newcommand{\rarrow}{\:\raisebox{0pt}{\tikz{\draw[->,solid,line width = 0.5 pt](0,0) -- (5mm,0);}}\:}

\newcommand{\nrarrow}{\:\raisebox{0pt}{\tikz{\draw[->,solid,line width = 0.5 pt](0,0) -- (5mm,0);
            \draw[-,solid,line width = 0.5 pt](1.5mm,-1mm) -- (2.5mm,1mm);}}\:}
\newcommand{\nrarrowzwei}{\:\raisebox{0pt}{\tikz{\draw[->,solid,line width = 0.5 pt](0,0) -- (5mm,0);
            \draw[-,solid,line width = 0.5 pt](1.5mm,-1mm) -- (2.5mm,1mm);}}}
\newcommand{\nrarrownull}{\nrarrowzwei_{0\:}}   
   
\newcommand{\nsimnull}{\nsim_{0\:}}

\allowdisplaybreaks 
\definecolor{darkgreen}{RGB}{0,139,0}
\newcommand{\LS}[1]{{\color{purple} #1}}
\newcommand{\VF}[1]{{\color{darkgreen} #1}}

\begin{document}

\begin{frontmatter}
\title{Mixed orthogonality graphs for continuous-time \vspace*{0.2cm} \\ state space models and orthogonal projections \vspace*{0.2cm} }
\runtitle{orthogonality graphs for ICCSS processes}

\begin{aug}
{   \author{\fnms{Vicky} \snm{Fasen-Hartmann}\ead[label=e1]{vicky.fasen@kit.edu}\orcid{0000-0002-5758-1999}}
   \and
    \author{\fnms{Lea} \snm{Schenk}\ead[label=e2]{lea.schenk@kit.edu}\ead[label=e3]{}\orcid{0009-0009-6682-6597}}
}
 \address{Institute of Stochastics, Karlsruhe Institute of Technology\\[2mm] \printead[presep={\ }]{e1,e2}}


\thankstext{e3}{\textsl{Funding:} This work is supported by the project “digiMINT”, which is a part of the “Qualitätsoffensive Lehrerbildung”, a joint initiative of the Federal Government and the Länder which aims to improve the quality of teacher training. The program is funded by the Federal Ministry of Education and Research. The authors are responsible for the content of this publication.}

\runauthor{V. Fasen-Hartmann and L. Schenk}
\end{aug}

\begin{abstract}
In this paper, we derive (local) orthogonality graphs for the popular continuous-time state space models,  including in particular multivariate continuous-time ARMA (MCARMA) processes.
In these (local) orthogonality graphs, vertices represent the components of the process, directed edges between the vertices indicate causal influences and undirected edges indicate contemporaneous correlations between the component processes. 
We present sufficient criteria for state space models to satisfy the assumptions of  \cite{VF23pre} so that the (local) orthogonality graphs are well-defined and various 
Markov properties hold. Both directed and undirected edges in these graphs are characterised by orthogonal projections on well-defined linear spaces. To compute these orthogonal projections, we use the unique controller canonical form of a state space model, which exists under mild assumptions, to recover the input process from the output process. We are then able to derive some alternative representations of the output process and its highest derivative. Finally, we apply these representations to calculate the necessary orthogonal projections, which culminate in the characterisations of the edges in the  (local) orthogonality graph. These characterisations 
are given by the parameters of the controller canonical form and the covariance matrix of the driving Lévy process.

\end{abstract}



\end{frontmatter}

\section{Introduction}\label{sec:intro} 

State space models are important tools in many scientific and engineering disciplines, including control theory, statistics, and computational neuroscience. In this paper, we study 
the time-invariant $\R^k$-valued \textsl{state space model} $(\Astar, \Bstar, \Cstar, L)$ of dimension $kp$, that is characterised by a driving $\R^k$-valued Lévy process $L=(L(t))_{t\in\R }$, a state transition matrix $\Astar \in \R^{kp\times kp}$ with $p\in\N$, an input matrix $\Bstar \in \R^{kp\times k}$, and an observation matrix $\Cstar \in \R^{kp\times k}$. Note that an $\R^k$-valued Lévy process $L$ is a stochastic process with stationary and independent increments, it is continuous in probability, and satisfies $L(0)=0_k\in \R^k$ almost surely \citep{SA07}.  A continuous-time state space model  then consists  of a state equation
\begin{align} \label{controller dgl}
dX(t)= \Astar X(t)dt+ \Bstar dL(t),
\end{align}
and an observation equation
\begin{align*}
   Y(t) = \Cstar X(t).
\end{align*}
The $\R^{kp}$-valued process $\CX=(X(t))_{t\in\R}$ is the input process and the $\R^{k}$-valued process $\CY=(Y(t))_{t\in\R}$ is the output process. 
It is well known that the solution of the state equation \eqref{controller dgl} is
\begin{align} \label{X aufgeteilt}
    X(t)=e^{\Astar (t-s)}X(s)+ \int_s^t e^{\Astar (t-u)} \Bstar dL(u), \quad s<t.
\end{align}
A special subclass of such state space models is the popular \textsl{multivariate continuous-time ARMA} (MCARMA) models \citep{MA07, SC12, SC122}. In contrast to \cite{SC12}, we speak here of a subclass instead of the equivalence of these classes because in our opinion there is an argument in the proof that is not clearly verifiable; the details are given in \Cref{subsec: state space models}.



This paper aims to construct a graphical model for such state space models.  The interest in graphical models for stochastic processes has increased significantly in recent years, see, e.g.,~\cite{Mogensen:Hansen:2022, Mogensen:Hansen:2020, Basu:Shojaie, EI07, EI10, DI07, DI08,VF23pre, VF23prec}, although the use of graphical models to visualise and analyse dependence structures in stochastic models is quite old \citep{Wright1921, Wright1934}. A major reason for this surge in interest is the simplicity and clarity of graphical models in representing the dependence structure in stochastic models such that examples of practical applications are ubiquitous. Another big advantage is their ease of implementation on computers, making them a powerful tool for the analysis of high-dimensional time series, as demonstrated, for example, in \cite{EI07}. The state of the art of graphical models is presented in \cite{Handbook:graphical}. 
 
In this paper, we use the approach of \cite{VF23pre} to construct \textsl{orthogonality graphs} and \textsl{local orthogonality graphs} for state space models, and to derive analytic representations of the edges in these graphs by the model parameters. 
(Local) orthogonality graphs are mixed graphs where the vertices $V=\{1,\ldots,k\}$ represent the different component series $\CY_v=(Y_v(t))_{t\in\R}$, $v\in V$, of an $\R^k$-valued process $Y_V=\CY$.  Directed edges reflect  (local) Granger causality and undirected edges reflect (local) contemporaneous correlation between the component series of the stochastic process. An attractive property of the (local) orthogonality graph is that it satisfies several types of Markov properties under fairly general assumptions. 

The mathematical notion of causality can be traced back to \cite{GR69} and \cite{SI72} and has since then been extended and applied in various fields, see \cite{Shojaie:Fox} for an excellent survey. In our context of continuous-time stochastic processes, the notion of Granger causality and contemporaneous correlation, as defined in \cite{VF23pre}, are based on conditional orthogonality relations of linear subspaces generated by subprocesses, similarly to \cite{EI07} in discrete time. This setup is perfectly suitable for stationary stochastic processes. At the same time, the approaches of \cite{EI10} using conditional independence relations for stochastic processes in discrete time and that of~\cite {Mogensen:Hansen:2022, Mogensen:Hansen:2020, DI07, DI08, Eichler:Dahlhaus:Dueck} using conditional local independence are suitable for semimartingales and point processes. We refer to \cite{VF23pre} for a detailed overview of graphical models for stochastic processes and the advantages of the different approaches.

The conditional orthogonality relations in the definition of (local) Granger causality and (local) contemporaneous correlation can be expressed equivalently by orthogonal projections of the component processes $Y_v(t+h)$ ($t\in\R$, $h\geq 0$, $v\in V$) and their highest derivative, respectively, on well-defined linear subspaces. To the best of our knowledge, the orthogonal projections for multivariate state space models and their derivatives required in this paper 
have not yet been addressed in the existing literature. Although \cite{RO67}, III, 5, is devoted to the topic of predictions for general stationary processes, 
the representations in that book are based on a specific maximal decomposition of the spectral density matrix of that process. This decomposition is generally not expressible as a simple formula, and so he only considers univariate examples. The orthogonal projections of univariate CARMA processes were discussed in the previous paper by \cite{BR15}. They provide representations for the linear projection of a CARMA process $Y(h)$ onto the entire linear space generated by the CARMA process up to time $t=0$,  and for the conditional expectation of $Y(h)$ on the $\sigma$-algebra generated by the CARMA process up to time $t=0$. A multivariate generalisation of the conditional expectation result using the $\sigma$-algebra generated by the whole multivariate CARMA (MCARMA) process up to time $t=0$ can be found in \cite{BA19}, Corollary 4.11, but the statement is not consistent with the comprehensible result of \cite{BR15}. In any case, in this paper, we require not only projections of $Y_v(t+h)$  on the linear space of the past of the process $\CY$ up to time $t$, but also on linear subspaces generated by subprocesses, and additionally the projections of the highest derivatives of $Y_v(t+h)$. 

In the context of \textsl{multivariate continuous-time AR} (MCAR) processes, which are a subclass of MCARMA processes,  the topic of orthogonal projections and the corresponding (local) orthogonality graphs are discussed in \cite{VF23pre}. Although MCAR processes are state space models, the techniques used there are not applicable to  MCARMA$(p,q)$ processes with $q>0$, because MCAR models have a much simpler structure. This structure allows,  e.g., the direct recovery of the input process $\CX$ from the output process $\CY$. In particular, \cite{VF23pre}  use the orthogonal projections to develop  (local) orthogonality graphs for MCAR processes and give interpretations for the edges which correspond to other results in the literature; see, e.g.,  \cite{CO96} for causality relations of MCAR processes and \cite{EI07} for discrete-time AR processes. The present paper can therefore be seen as an extension of \cite{VF23pre} to a broader class of models and of course, we compare our results with those of that paper. We would like to point out here that, according to our knowledge, even for stochastic processes in discrete time, the literature on mixed graphical models is restricted to AR processes \citep{EI07, EI10}, there is not much known on mixed graphical models for the more complex ARMA processes satisfying some types of Markov properties.

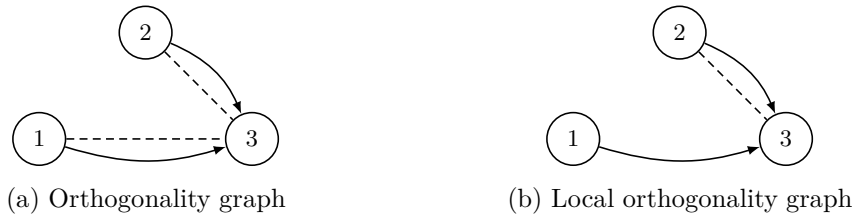
\begin{figure}[ht]
  \begin{subfigure}{0.49\textwidth}
  \centering
    \begin{tikzpicture}[align=center,node distance=2cm and 2cm, semithick , state/.style ={circle, draw,  text=black , minimum width =0.7 cm}, arrow/.style={-latex}, every loop/.style={}]

  		\node[state] (2) {2};
  		\node[state] (1) [below left of=2] {1};
  		\node[state] (3) [below right of=2] {3};

  		\path
   		(1) edge [densely dashed] node {} (3)
   		(2) edge [densely dashed] node {} (3)
            
   		(1) edge [arrow, bend right=17] node {} (3)
            (2) edge [arrow, bend left=22] node {} (3)
		;
 	\end{tikzpicture}  
  \caption{Orthogonality graph}
  \label{fig: Orthogonality graph a}
  \end{subfigure}
  \begin{subfigure}{0.49\textwidth}
  \centering
	\begin{tikzpicture}[align=center,node distance=2cm and 2cm, semithick ,
		state/.style ={circle, draw,  text=black , minimum width =0.7 cm}, arrow/.style={-latex}, every loop/.style={}]

  		\node[state] (2) {2};
  		\node[state] (1) [below left of=2] {1};
  		\node[state] (3) [below right of=2] {3};

  		\path
   		(2) edge [densely dashed] node {} (3)
     
   		(1) edge [arrow, bend right=17] node {} (3)
   		(2) edge [arrow, bend left=22] node {} (3)
		;
 	\end{tikzpicture}  
  \caption{Local orthogonality graph}
  \label{fig: Partial correlation graph b}
  \end{subfigure}
  \caption{In the left figure is the orthogonality graph and in right figure the local orthogonality graph of the $3$-dimensional ICCSS$(2,1)$ process defined in \Cref{example}.}
  \label{fig: Two graphical models b}
\end{figure}

In the present paper, the controller canonical form of a state space model plays an important role in calculating the orthogonal projections of $Y_v(t+h)$ and its highest derivative, as highlighted in \cite{BA19}.  The controller canonical form of an MCARMA process has been studied in \cite{BR12} and is the multivariate generalisation of the definition of a univariate CARMA process \citep{BR14}. In the case of the existence of a controller canonical form, we show that this representation is unique, which is essential for the unique characterisations of the edges in the (local) orthogonality graph later in this paper; we will prove that the edges depend only on the model parameters of the controller canonical form and the covariance matrix of the driving Lévy process. A special feature of the controller canonical form is that under very general assumptions it can be used to recover the input process $\CX$ from the output process $\CY$, i.e., $X(t)$ lies in the closed linear space generated by $(Y(s))_{s\leq t}$ and which has only been known for univariate CARMA processes \cite[Theorem 2.2]{BR15}. In this case, we call $\CY$ an \textsl{invertible controller canonical state space} (ICCSS) process and we are able to calculate the necessary orthogonal projections to describe the  (local) Granger causality and (local) contemporaneous correlation relations of the underlying process, resulting in the characterisations of the edges in the (local) orthogonality graph. 
An example of both graphs for a 3-dimensional ICCSS$(2,1)$ process is given in \Cref{fig: Two graphical models b}. This example is used for illustration throughout the paper. 

In conclusion, we show in this paper that not only for ICCSS processes but also for most state space models, the (local) orthogonality graph exists and satisfies several of the preferred Markov properties. In our opinion, this is the first (mixed) graphical model for this popular and broad class of stochastic processes. However, for the explicit representation of the edges via (local) Granger causality and (local) contemporaneous correlation,  we need the invertibility of the state space model and hence, the restriction to ICCSS processes. In addition, we derive new and important results for state space models, such as sufficient criteria for being an ICCSS model, alternative representations for $Y_v(t+h)$ and its highest derivative depending on the linear past of $\CY$ up to time $t$ and an independent noise term, and, in particular, their orthogonal projections onto linear subspaces, which are the basis for linear predictions.

The paper is structured as follows. In \Cref{sec: Preliminaries}, we introduce (local) orthogonality graphs, the controller canonical form of a state space model and ICCSS processes. We also present their basic properties that are important for this paper. In \Cref{sec: Linear prediction of MCARMA processes}, we consider orthogonal projections of ICCSS processes and their highest derivatives onto linear subspaces generated by subprocesses, since the characterisations of the edges of the (local) orthogonality graph are based on these orthogonal projections. These results then lead to the existence of the (local) orthogonality graph and the characterisation of the directed and the undirected edges of an ICCSS process by its model parameters in \Cref{sec: orthogonality graph for MCARMA processes}.  The proofs of the paper are given in their own \Cref{Appendix:Proofs}.

\subsection*{Notation}
From now on we call the space of all real and complex $(k\times k)$-dimensional matrices $M_k(\R)$ and $M_k(\C)$, respectively. Similarly, $M_{k,d}(\R)$ and $M_{k,d}(\C)$ denote real and complex $(k\times d)$-dimensional matrices. For $\BA \in M_{k,d}(\C)$ we write $\BA^\top$ for the transpose of $\BA$ and
for $\BA \in M_{k}(\C)$ we write
$\BA\geq 0$ if $\BA$ is positive semi-definite, and $\BA>0$ if $\BA$ is positive definite. Furthermore, $\sigma(\BA)$ are the eigenvalues of $\BA$. $I_{k}$ is the $(k\times k)$-dimensional identity matrix, $0_{k\times d}$ is the $(k\times d)$-dimensional zero matrix and $0_{k}$ is either the $k$-dimensional zero vector or the $(k\times k)$-dimensional zero matrix which should be clear from the context. The vector $e_a \in \R^k$ is the $a$-th unit vector,  as well as
\begin{gather}\label{eq: def Ej}
\BFE_j  =
\begin{pmatrix}
     0_{k(j-1) \times k} \\
     I_{k} \\
     0_{k(p-j) \times k}
\end{pmatrix} \in M_{kp \times k}(\R), \quad j=1,\ldots,p.
\end{gather} 
Furthermore, we write for $p>q$,
\begin{gather}\label{eq: def E}
\OBE^\top = (I_{kq} , 0_{k(p-q)\times kq}) \in M_{kq \times kp}(\R) \quad \text{ and } \quad
\UBE^\top = (0_{k \times k(q-1)}, I_k) \in M_{k \times kq}(\R).
\end{gather}
For a matrix polynomial $P(z)$ we define the set of zeros of the polynomial $\det (P(z))$ as $\mathcal{N}(P)=\{z\in \C:\det (P(z))=0\}$ and $\text{deg}(\det(P(z)))$ denotes the degree of the polynomial $\det(P(z))$. 
$\CY=(Y(t))_{t\in \R}$ is a $k$-dimensional stationary and mean-square continuous stochastic process with index set $V=\{1,\ldots, k\}$ and expectation zero. The corresponding components are denoted by $\CY_v=(\CY_v(t))_{t\in \R}$ for $v\in V$ and  multivariate subprocesses are denoted by $\CY_S = (\CY_s)_{s\in S} = (Y_S(t))_{t\in \R}$ for $S\subseteq V$. A special case is $S=V$ where $Y_V=Y$. Finally, $\text{l.i.m.}$ is the mean square limit. 

\section{Preliminaries}\label{sec: Preliminaries}

The topic of this paper is mixed (local) orthogonality graphs for state space models. Therefore, in \Cref{subsec: orthogonality graphs} we introduce (local) orthogonality graphs and present the main results of \cite{VF23pre} that are relevant to this paper. Then, in Sections \ref{subsec: state space models} and \ref{subsec: Inversion of MCARMA processes},  the controller canonical form of a state space model is introduced, and the invertibility of these processes is discussed. 

\subsection{Mixed orthogonality graphs}\label{subsec: orthogonality graphs}

Orthogonality graphs are graphical models \linebreak  $G_{OG}=(V,E_{OG})$ where the vertices $V=\{1,\ldots,k\}$ represent the different component series  $\CY_v = (Y_v(t))_{t\in \R}$, $v \in V$, of a $k$-dimensional wide-sense stationary and mean-square continuous stochastic process $\CY=(Y(t))_{t\in \R}$ with expectation zero. The vertices are connected with both directed and undirected edges $E_{OG}$. A directed edge represents a Granger causal relationship while an undirected edge represents a contemporaneous correlation between the components. To make this clear, we first present some notations that we require for the definitions of Granger causality and contemporaneous correlation.  For any set $S\subseteq V$, $t, \widetilde{t} \in \R$ with $t<\widetilde{t}$, we define the closed linear spaces generated by the  subprocess $\CY_S =(Y_s)_{s\in S}= (Y_S(t))_{t\in \R}$ as
\begin{align*}
    \mathcal{L}_{Y_S}(t, \widetilde{t}) &= \overline{\Bigg\{ \sum_{i=1}^n \sum_{s\in S} \gamma_{s,i} Y_s(t_i): \gamma_{s,i} \in \C,\: t \leq t_1 \leq \cdots \leq t_n \leq \widetilde{t}, n \in \N  \Bigg\}},\\
    \mathcal{L}_{Y_S}(t) &= \overline{\Bigg\{ \sum_{i=1}^n \sum_{s\in S} \gamma_{s,i} Y_s(t_i): \gamma_{s,i} \in \C, \: -\infty < t_1 \leq \cdots \leq t_n \leq t, n \in \N  \Bigg\}}.
\end{align*}
We further denote the orthogonal projection of $Z\in L^2$ on the closed linear subspace $\mathcal{L}\subseteq L^2$ by $P_{\mathcal{L}}(Z)=P_{\mathcal{L}}Z$.
Now, we establish definitions of  Granger causality, which characterises the directed edges in the (local) orthogonality graph.

\begin{definition}\label{  Granger non-causal} \label{Charakterisierung als Gleichheit der linearen Vorhersage}
Let $a,b\in S \subseteq V$ and $a \neq b$. 
\begin{itemize}
    \item[(a)] $\CY_a$ is \textsl{Granger non-causal for $\CY_b$ with respect to $\CY_S$}, if and only if,  
    for all $t\in \R$ and $0 \leq h \leq 1$,
\begin{align*}
    &P_{\mathcal{L}_{Y_S}(t)}Y_b(t+h) = P_{\mathcal{L}_{Y_{S \setminus \{a\}}}(t)}Y_b(t+h) \quad \text{$\mathbb{P}$-a.s.~} \nonumber 
\end{align*}
We shortly write $\CY_a \nrarrow \CY_b \: \vert \: \CY_S$.
\item[(b)] Suppose  $\CY_v$ is $j_v$-times mean-square differentiable
but the $(j_v+1)$-derivative does not exist for  $v\in V$. The $j_v$-derivative is denoted by $D^{(j_v)}\CY_v$, where  for $j_v=0$ we define $D^{(0)} \CY_v=\CY_v$. Then $\CY_a$ is \textsl{locally Granger non-causal} for $\CY_b$ with respect to $\CY_S$, if and only if, for all $t\in \R$,

\begin{align*}
&\limh P_{\mathcal{L}_{Y_S}(t)} \left(\frac{D^{(j_b)} Y_b(t+h)- D^{(j_b)} Y_b(t)}{h}\right)\\
&\quad =\limh P_{\mathcal{L}_{Y_{S \setminus \{a\}}}(t)} \left(\frac{D^{(j_b)} Y_b(t+h)- D^{(j_b)} Y_b(t)}{h}\right)
 \quad \mathbb{P}\text{-a.s.}
\end{align*}
 We shortly write $\CY_a \nrarrownull \CY_b \: \vert \: \CY_S$.
 \end{itemize}
\end{definition}

In words,  $\CY_a$ is Granger non-causal for $\CY_b$ with respect to $\CY_S$ if the prediction of $Y_b(t+h)$ based on the linear information available at time $t$ provided by the past and present of $\CY_S$ is not diminished by removing the linear information provided by the past and present values of $\CY_a$. Local Granger non-causality considers the limiting case $h\rightarrow 0$, where the highest existing derivative of the process must be examined to obtain a non-trivial criterion; see as well the discussion in \cite{VF23pre} for the motivation of these definitions. 

 \begin{remark}  
 \cite{VF23pre} originally defined Granger causality by conditional orthogonality of linear spaces generated by subprocess, and then showed that these definitions are equivalent to the definitions based on the orthogonal projections given above (see Theorem 3.5 and Remark 3.12 therein). 
\end{remark}
Next, the undirected edges are characterised by (linear) contemporaneous correlation. The idea is simple: there is no undirected influence between $\CY_a$ and $\CY_b$ with respect to $\CY_S$ if and only if, given the amount of information provided by the past of $\CY_S$ up to time $t$, $\CY_a$ and $\CY_b$ are uncorrelated in the future. \cite{VF23pre} introduce the following definitions.

\begin{definition} \label{Charakterisierung als Gleichheit der linearen Vorhersage 2} \label{Def: local cont uncorr}
Let $a,b\in S \subseteq V$ and $a \neq b$. 
\begin{itemize}
    \item[(a)] $\CY_a$ and $\CY_b$ are \textsl{contemporaneously uncorrelated} with respect to $\CY_S$, if and only if, for all $t\in \R$ and $0\leq h, \widetilde{h} \leq  1$,
\begin{align*}
\BE & \left[\left( Y_a(t+h)- P_{\mathcal{L}_{Y_{S}}(t)} Y_a(t+h) \right)\overline{\left( Y_b(t+\widetilde{h})- P_{\mathcal{L}_{Y_{S}}(t)} Y_b(t+\widetilde{h}) \right)} \right]=0.
\end{align*}
We shortly write $\CY_a \nsim \CY_b \: \vert \: \CY_S$.
\item[(b)] Suppose  $\CY_v$ is $j_v$-times mean-square differentiable but the $(j_v+1)$-derivative does not exist for  $v\in V$. 
Then $\CY_a$ and $\CY_b$ are \textsl{locally contemporaneously uncorrelated} with respect to $\CY_S$, if and only if, for all $t\in \R$, 
\begin{align*}
\limhh \frac{1}{h} \: \BE & \bigg[\left( D^{(j_a)} Y_a(t+h)- P_{\mathcal{L}_{Y_{S}}(t)} D^{(j_a)} Y_a(t+h) \right)  \\
		  &\times \overline{\left( D^{(j_b)} Y_b(t+h)- P_{\mathcal{L}_{Y_{S}}(t)} D^{(j_b)} Y_b(t+h) \right)} \bigg]=0.
\end{align*}
 We shortly write $\CY_a \nsimnull \CY_b \: \vert \: \CY_S$.
\end{itemize}
\end{definition}

\begin{remark}
The original definition of contemporaneous uncorrelatedness in \cite{VF23pre} is again by conditional orthogonality of linear spaces generated by subprocesses, and equivalent characterisations using orthogonal projections are given there (see Lemma 4.3, Theorem 4.5, and Remark 4.9 therein).
\end{remark}

Before defining the orthogonality graphs with these terms and definitions, we introduce assumptions on the stochastic process $\CY$, which are fulfilled, in particular, by most state space models (see \Cref{sec: orthogonality graph for MCARMA processes}).

\pagebreak

\begin{assumptionletter} \label{Assumption graph} $\mbox{}$
\renewcommand{\theenumi}{(\ref*{Assumption graph}.\arabic{enumi})}
         \renewcommand{\labelenumi}{\theenumi}
\begin{enumerate}
    \item \label{stationarity}
        $\CY$ is a $ k$-dimensional wide-sense stationary and mean-square continuous stochastic process with expectation zero and index set $V=\{1,\ldots,k\}$. 
    \item \label{Assumption an Dichte}
        $\CY$ has a spectral density matrix $f_{YY}(\lambda)>0$ and   there exists an $0<\varepsilon<1$, such that 
        \begin{align*}
            (1-\varepsilon)I_{\alpha} - f_{Y_AY_A}(\lambda)^{-1/2}f_{Y_AY_B}(\lambda)f_{Y_BY_B}(\lambda)^{-1}f_{Y_BY_A}(\lambda)f_{Y_AY_A}(\lambda)^{-1/2} \geq 0,
        \end{align*}
        for almost all $\lambda \in \R$ and for all disjoint subsets $A,B\subseteq V$, where $\alpha$ is the cardinality of $A$.
    \item \label{Assumption purely nondeterministic of full rank} 
        $\CY$ is purely non-deterministic, i.e., for all $a\in V$ and $t\in \R$,
        \begin{align*}
            \underset{h \rightarrow \infty}{\text{l.i.m.\:}} P_{\mathcal{L}_{Y}(t)} Y_a(t+h) = 0 \quad \text{$\mathbb{P}$-a.s.}
        \end{align*}
    \end{enumerate}
\end{assumptionletter}

\begin{remark} $\mbox{}$
\begin{itemize} 
    \item[(a)] Assumption \ref{stationarity} is a basic requirement, otherwise, e.g., the spectral density in Assumption \ref{Assumption an Dichte} is not well defined.  
    \item[(b)] Assumption \ref{Assumption an Dichte} ensures the linear independence and closedness of sums of linear spaces generated by subprocesses, i.e.,~for $t\in \R$ and disjoint subsets $A, B\subseteq V$,
\begin{align}\label{eq: linear independence}
\mathcal{L}_{Y_{A}}(t) \cap \mathcal{L}_{Y_{B }}(t) =  \{0\} \quad \text{and} \quad
\mathcal{L}_{Y_{A}}(t) + \mathcal{L}_{Y_{B}}(t) = \overline{\mathcal{L}_{Y_{A}}(t) + \mathcal{L}_{Y_{B}}(t)} \quad \mathbb{P}\text{-a.s.}
\end{align}
    \item[(c)] Any process that is wide-sense stationary can be uniquely decomposed into a deterministic and a purely non-deterministic process, which are orthogonal to each other (\citealp{GL58}, Theorem 1). From the point of view of applications, deterministic processes are not important, so \ref{Assumption purely nondeterministic of full rank} is a natural assumption.
\end{itemize}
\end{remark}

With this assumption, we are able to define orthogonality graphs.

\begin{definition}\label{Definition orthogonality graph}
Suppose $\CY$ 
satisfies \Cref{Assumption graph}. 
\begin{itemize}
    \item[(a)]
If we define $V=\{1,\ldots,k\}$ as the vertices and the edges $E_{OG}$ via
\begin{itemize}
\item[(i)\phantom{i}]   $a \rarrow b \notin E_{OG} \quad \Leftrightarrow \quad \CY_a \nrarrow \CY_b \: \vert \: \CY_V$,
\item[(ii)]  $a \inst b \notin E_{OG} \quad \Leftrightarrow \quad \CY_a \nsim \CY_b \: \vert \: \CY_V$,
\end{itemize}
for $a,b\in V$ with $a\neq b$, then $G_{OG}=(V,E_{OG})$ is called \textsl{orthogonality graph} for $Y=\CY_V$.
\item[(b)] If we define $V=\{1,\ldots,k\}$ as the vertices and the edges $E_{OG}^{0}$ via
\begin{itemize}
\item[(i)\phantom{i}]   $a \rarrow b \notin E_{OG}^{0} \quad \Leftrightarrow \quad \CY_a \nrarrownull \CY_b \: \vert \: \CY_V$,
\item[(ii)]  $a \inst b \notin E_{OG}^{0} \quad \Leftrightarrow \quad \CY_a \nsimnull \CY_b \: \vert \: \CY_V$,
\end{itemize}
for $a,b\in V$ with $a\neq b$, then $G_{OG}^{0}=(V,E_{OG}^{0})$ is called \textsl{local orthogonality graph} for $Y=\CY_V$.
\end{itemize}
\end{definition}

\begin{remark}\label{rem:Markov properties}
\begin{itemize} $\mbox{}$
   \item[(a)] As discussed in \cite{VF23pre}, the assumptions are not necessary for the definition of the graphs, but they ensure that the usual Markov properties for mixed graphs are satisfied. Specifically, the orthogonality graph satisfies the pairwise, local, block-recursive, global AMP and global Granger-causal Markov properties. In particular,  Assumption \ref{Assumption purely nondeterministic of full rank} ensures that the global AMP Markov property holds. The local orthogonality graph satisfies the pairwise, local and block-recursive Markov property; for global Markov properties of the local orthogonality graph additional assumptions are required; see \cite{VF23pre}, Propositions 5.20 and 5.21, respectively.
   \item[(b)] The local orthogonality graph has fewer edges than the orthogonality graph and, in general, the graphs are not equal. An explicit example of a (local) orthogonality graph, illustrating this property, is given in \Cref{fig: Two graphical models b}. The advantage of the local orthogonality graph over the orthogonality graph is that it allows for modeling more general graphs, whereas the orthogonality graph satisfies the global AMP and the global causal Markov property, the local orthogonality graph does not satisfy them in general, additional assumptions are necessary. For more details see \cite{VF23pre} again.
\end{itemize}
\end{remark}

\subsection{Controller canonical state space models}\label{subsec: state space models}
The paper aims to derive (local) orthogonality graphs for state space models. Therefore, we use the controller canonical form of a state space model and its uniqueness, which results in the unique characterisation of the edges in the (local) orthogonality graph in \Cref{sec: orthogonality graph for MCARMA processes}. To define the controller canonical form of a state space model, we need some definitions and terminology. Therefore, note that each state space model $(\Astar, \Bstar, \Cstar, L)$ is associated with a rational matrix function 
\begin{align}\label{eq: Transfer function}
    H(z)=\Cstar \left(z I_{kp} - \Astar \right)^{-1} \Bstar,
    \quad z\in\C\backslash \sigma(\Astar),
\end{align}
called the \textsl{transfer function} of the state space model, and the triple $(\Astar, \Bstar, \Cstar)$ is an \textsl{algebraic realisation} of the transfer function of dimension $kp$ (characterising the dimension $(kp\times kp)$ of $\Astar$). The triple $(\Astar, \Bstar, \Cstar)$ is said to be \textsl{minimal} if there is no other algebraic realisation of $H(z)$ with dimension less than $kp$. The transfer function is of importance because, due to the spectral representation theorem \cite[Theorem 17.5]{lax2002}, we are able to recover the kernel function $\Cstar\mbox{e}^{\Astar t}\Bstar{\boldsymbol 1}_{\{t\geq 0\}}$ of the output process $Y$ of the state space model via
\begin{align}
    \Cstar\mbox{e}^{\Astar t}\Bstar=\frac{1}{2\pi i}\int_\Gamma \mbox{e}^{zt} H(z)\,dz, \quad t\geq 0,
\end{align}
where $\Gamma$ is a closed contour in the complex numbers that winds around each eigenvalue of $\Astar$ exactly once. The transfer function even uniquely determines the function $\Cstar\mbox{e}^{\Astar t}\Bstar$, $t\in\R$ \citep[Lemma~3.2]{SC122}.
 \cite{KA80} provides in Lemma 6.3-8 that there exist $(k\times k)$-dimensional matrix polynomials $P(z)$ and $Q(z)$ such that
\begin{align}\label{eq: coprime right polynomial fraction description}
    H(z) = Q(z)P(z)^{-1}, \quad z\in \C \setminus \sigma(\Astar) ,
\end{align}
is a \textsl{coprime} right polynomial fraction description of the transfer function, 
which in turn means that the matrix $[P(z) \:\: Q(z)]$ has full rank for all $z\in\C$. 
In Lemma 6.3-3 \cite{KA80} even gives a construction for such a decomposition. However, without any additional assumption, the coprime polynomials $P(z)$ and $Q(z)$ that satisfy \eqref{eq: coprime right polynomial fraction description} are not unique. Even if we assume additionally that $(\Astar, \Bstar, \Cstar)$ is minimal, resulting in $\deg(\det(P(z))=kp$  \cite[Theorem 17.5]{RU96}, then we can take any invertible matrix $S\in M_{k}(\R)$ such that the matrix polynomials $P(z) S$ and $Q(z) S$ also satisfy $H(z)= Q(z) S S^{-1} P(z)^{-1}$.


Despite the many different coprime polynomials $P(z)$ and $Q(z)$ that satisfy \eqref{eq: coprime right polynomial fraction description}, to the best of our knowledge, it remains unclear whether there is a coprime right polynomial fraction description with 
\begin{align}\label{def PQ}
P(z)=I_k z^p+A_1z^{p-1}+\ldots+A_p \quad \text{ and } \quad
Q(z)=C_0+C_1 z+\ldots+C_q z^q, 
\end{align}
$A_1 , A_2,\ldots, A_p, C_0,C_1,\ldots, C_q \in M_k(\R)$, and $p,q\in \N_0$, $p>q$, i.e., $z^p$ is the highest power with prefactor $I_k$. 
In  representation \eqref{def PQ}, the assumption  $\deg(\det(P(z))=kp$ obviously holds.
Note, the construction method of \cite{KA80} often gives a polynomial $P(z)$ with higher order than $p$, but the prefactor of the highest power has a zero determinant. \cite{BR12}, Theorem 3.2, and \cite{SC12}, Corollary~3.4, implicitly assume such a right polynomial fraction description \eqref{def PQ} without discussing its existence. Since the existence of such a coprime right polynomial fraction description is essential for the forthcoming results, we always assume it additionally. In Examples \ref{example} and \ref{example coprime} we present likewise examples where this assumption is fulfilled. 

\begin{example}\label{example} 
Consider a state space model $(\Astar, \Bstar, \Cstar, L)$ with $k=3$, $p=2$, $q=1$ and $\Sigma_L=I_3$,  where we set
\setlength{\arraycolsep}{4pt}
\begin{align*}
    \Astar=\begin{pmatrix*}[r]
        0 & -1 & 0 & 0 & 0 & 0\\
        1 & -1 & 0 & 0 & 0 & 0\\
        0 & 0 & 0 & -1 & 0 & 0\\
        0 & 0 & 1 & -1 & 0 & 0\\
        0 & 0 & 0 & - \frac{1}{2} & 0 & -1\\
        0 & 0 & 0 & 0 & 1 & -1
    \end{pmatrix*}, \quad
    \Bstar=\begin{pmatrix*}[r]
        1 & -1 & 0\\
        1 & -1 & 0\\
        0 & 2 & 0\\
        0 & 2 & 0\\
        0 & -1 & 1\\
        0 & -1 & 1
    \end{pmatrix*} \quad \text{ and }\quad
    \Cstar=\begin{pmatrix}
        0 & 1 & 0 & 0 & 0 & 0 \\
        0 & 1 & 0 & 1 & 0 & 0 \\
        0 & 1 & 0 & \frac{3}{2} & 0 & 1
    \end{pmatrix}.
\end{align*}
Note that this is the state space model that generates the (local) orthogonality graph in \Cref{fig: Two graphical models b} in the introduction to this paper, and we will look at this example in more detail in the course of this paper.

In this example, a straightforward calculation shows that there exists the right polynomial fraction description of the transfer function
\setlength{\arraycolsep}{4pt}
\begin{align*}
    H(z)
    &= \Cstar \left(z I_{6} - \Astar \right)^{-1} \Bstar\\
    &= \begin{pmatrix}
        z + 1 & -z - 1 & 0 \\
        z + 1 & z + 1 & 0 \\
        z + 1 & z + 1 & z + 1
    \end{pmatrix}
    \begin{pmatrix}
         z^2 + z + 1 & 0 & 0 \\
         0 & z^2 + z + 1 & 0 \\
         0 & 1 & z^2 + z + 1
    \end{pmatrix}^{-1}\\
    &=: Q(z) P(z)^{-1},
\end{align*}
for $z\in \C \setminus \sigma(\Astar)$, where $\sigma(\Astar)=\{-\frac{1}{2}+ \frac{\sqrt{3}}{2}i, -\frac{1}{2}- \frac{\sqrt{3}}{2}i\}$. Furthermore, this decomposition is coprime, since $[P(z) \:\: Q(z)]$ has full rank 3 for all $z\in \C$. Thus, in this example, there exists a coprime right polynomial fraction description  \eqref{eq: coprime right polynomial fraction description} with polynomials as in \eqref{def PQ}. 
\end{example}

For the purpose of this paper, not only the existence of a coprime right polynomial fraction description of the transfer function with polynomials $P(z)$ and $Q(z)$ as in \eqref{def PQ} is important but also its uniqueness. In the next proposition, we derive that this requirement is immediately satisfied under the existence assumption of $P(z)$ and $Q(z)$. 



 
\begin{proposition}\label{lemma: uniqueness}
Let $(\Astar, \Bstar, \Cstar, L)$ be a 
state space model with transfer function $H(z)$. Suppose there exists a coprime right polynomial fraction description 
of $H(z)$ with polynomials $P(z)$ and $Q(z)$ as in \eqref{def PQ} such that
\begin{align*}
    H(z)=\Cstar \left(z I_{kp} - \Astar \right)^{-1} \Bstar = Q(z)P(z)^{-1}, \quad z\in \C \setminus \sigma(\Astar).
\end{align*}
Then $P(z)$ and $Q(z)$ are unique. Moreover, defining
\begin{align} \label{CCform}
\begin{split}
\BA &=
\begin{pmatrix*}[c]
  0_k & I_k & 0_k & \cdots & 0_k \\
  0_k & 0_k & I_k & \ddots & \vdots \\
  \vdots &  & \ddots & \ddots & 0_k \\
  0_k & \cdots & \cdots & 0_k & I_k \\
  -A_p & -A_{p-1} & \cdots & \cdots & -A_1
\end{pmatrix*}
\in M_{kp}(\R), \quad
\BB=
\begin{pmatrix}
0_k \\
\vdots \\
0_k \\
I_k
\end{pmatrix}
\in M_{kp\times k}(\R), \\
\BFC &=\begin{pmatrix}
C_0, & C_1, & \cdots, & C_{q}, & 0_k, & \cdots, & 0_k
\end{pmatrix} \in M_{k\times kp}(\R),
\end{split}
\end{align}
then $\sigma(\Astar)=\sigma(\BA)$ and
\begin{align*}
     H(z)=\BFC \left(z I_{kp} - \BA \right)^{-1} \BB, \quad z\in \C \setminus \sigma(\BA).
\end{align*}
Finally, $\CY$ is a solution of the state space model  $(\Astar, \Bstar, \Cstar, L)$, if and only if, it is a solution of the state space model  $(\BA, \BB, \BFC, L)$. The state space model $(\BA, \BB, \BFC, L)$ is called \textsl{controller canonical form}.
\end{proposition}

In particular, of course, this implies that there exists no other minimal state space representation with matrices of the structure as in \eqref{CCform}; this representation is unique. Since the solution of these two state space models is equal, we will henceforth assume, without loss of generality, that the state space model is given in the unique controller canonical form $(\BA, \BB, \BFC, L)$ as in \eqref{CCform}. 

\begin{example} \label{Example 2.10}
In \Cref{example}, the state space model $(\Astar, \Bstar, \Cstar, L)$ has the unique controller canonical form $(\BA, \BB, \BFC, L)$, where 
\setlength{\arraycolsep}{4pt}
\begin{align*}
    \BA = \begin{pmatrix*}[r]
        0 & 0 & 0 & 1 & 0 & 0 \\
        0 & 0 & 0 & 0 & 1 & 0 \\
        0 & 0 & 0 & 0 & 0 & 1 \\
        -1 & 0 & 0 & -1 & 0 & 0 \\
        0 & -1 & 0 & 0 & -1 & 0 \\
        0 & -1 & -1 & 0 & 0 & -1
        \end{pmatrix*}, \quad
    \BB = \begin{pmatrix}
        0 & 0 & 0 \\
        0 & 0 & 0 \\
        0 & 0 & 0 \\
        1 & 0 & 0 \\
        0 & 1 & 0 \\
        0 & 0 & 1
        \end{pmatrix}, \quad
     \BFC = \begin{pmatrix*}[r]
        1 & -1 & 0 & 1 & -1 & 0 \\
        1 & 1 & 0 & 1 & 1 & 0 \\
        1 & 1 & 1 & 1 & 1 & 1
        \end{pmatrix*}.
\end{align*}    
\end{example}

A particular example of a state space model is the MCARMA model whose output process $\CY$ is an MCARMA process. The definition of an MCARMA process $ \CY$  is motivated by the idea that
$ \CY$ solves the stochastic differential equation
\begin{align*}
    \Pstar(D) Y(t)=\Qstar(D)D L(t)
\end{align*}
where $D$ is the differential operator with respect to $t$ and
\begin{align} \label{PoQP}
 \Pstar(z) = I_k z^{\pstar} + \Pstar_1 z^{{\pstar}-1} + \ldots + \Pstar_{\pstar}  
   \quad \text{ and } \quad
  \Qstar(z) = \Qstar_0 z^{\qstar} +\Qstar_1 z^{{\qstar}-1} + \ldots + \Qstar_{\qstar}
\end{align}
are the AR (autoregressive) and MA (moving average) polynomial with $\Pstar_1, \Pstar_2,\ldots, \Pstar_{\pstar}$, $\Qstar_0, \Qstar_1,\ldots,\Qstar_{\qstar} \in M_k(\R)$. However, a Lévy process is not differentiable, so this is not a formal definition. The formal definition is given by a state space model \citep{MA07} as follows.

\begin{definition} \label{def:MCARMA}
Define
\begin{align*}
\Astar &= \begin{pmatrix}
  0_k & I_k & 0_k & \cdots & 0_k \\
  0_k & 0_k & I_k & \ddots & \vdots \\
  \vdots &  & \ddots & \ddots & 0_k \\
  0_k & \cdots & \cdots & 0_k & I_k \\
  -\Pstar_{\pstar} & -\Pstar_{\pstar-1} & \cdots & \cdots & -\Pstar_1
\end{pmatrix}
\in M_{k\pstar}(\R), \quad
\Bstar 
= \begin{pmatrix}
\beta_1 \\
\beta_2 \\
\vdots \\
\beta_{\pstar}
\end{pmatrix} 
\in M_{k\pstar\times k}(\R), \\
\Cstar &= (I_k , 0_k , \cdots , 0_k) \in M_{k\times k\pstar}(\R),
\end{align*}
where $\beta_{1}=\cdots=\beta_{\pstar-\qstar-1}= 0_k$ and $\beta_{\pstar-j} = -\sum_{i=1}^{\pstar-j-1} \Pstar_i \beta_{\pstar-j-i} + \Qstar_{\qstar-j}$ for $j=\qstar, \ldots,0$. Then $(\Astar, \Bstar, \Cstar, L)$ is called a \textsl{multivariate continuous-time moving average  model  of order $(\pstar,\qstar)$}, shortly an MCARMA$(\pstar,\qstar)$ model.
\end{definition}

\begin{remark}
A comparison of the triplets $(\BA, \BB, \BFC)$ and $(\Astar, \Bstar, \Cstar)$ shows that the MCAR$(p)=$ MCARMA$(p,0)$ model is already in controller canonical form. For MCARMA$(p,q)$ models, \cite{SC12} show the equivalence between state space models and MCARMA models in Corollary 3.4, and \cite{BR12}, Theorem 3.2, show the equivalence between MCARMA models and controller canonical state space models. However, as mentioned above, both implicitly assume the existence of a coprime left (right) polynomial fractional description \eqref{eq: coprime right polynomial fraction description} with polynomials $P(z)$ and $Q(z)$ as in \eqref{def PQ}, which is in our opinion not obvious. However, for univariate state space processes with $k=1$, the existence of a coprime right polynomial fractional description is apparent (see proof of \Cref{Xq als integral}), so that any univariate state space model is a CARMA model and vice versa; additionally, any univariate state space model has a representation in controller canonical form.
\end{remark}



A peculiarity of MCARMA models is that the AR polynomial $\Pstar(z)$ and the MA polynomial $\Qstar(z)$ provide a left polynomial fraction description of the transfer function, i.e., $H(z) = \Pstar(z)^{-1} \Qstar(z)$ (\citealp{MA07} or \citealp{BR12}, Lemma 3.1). If the MCARMA model is minimal, this left polynomial fractional description is even coprime \cite[Theorem 6.5-1]{KA80}. The connection to the coprime right polynomial fraction description \eqref{eq: coprime right polynomial fraction description} with $P(z)$ and $Q(z)$ as in \eqref{def PQ} is given in the next lemma.

\begin{lemma}\label{Lemma: Zusammenhang MCARMA}
Let an MCARMA$(\pstar,\qstar)$ model be given with state space representation
 $(\Astar, \Bstar, \Cstar, L)$ as in \Cref{def:MCARMA} and polynomials $\Pstar(z)$ and $\Qstar(z)$ as in \eqref{PoQP}.
Suppose $(\Astar, \Bstar, \Cstar)$ is minimal and there exists a coprime right polynomial fraction description \eqref{eq: coprime right polynomial fraction description} of the transfer function with polynomials $P(z)$ and $Q(z)$ as in \eqref{def PQ}.
Then $P(z)$ and $Q(z)$ are unique, $\pstar=p$, $\qstar=q$, $\Qstar_0=C_q$, $\mathcal{N}(\Pstar)=\mathcal{N}({P})$ and $\mathcal{N}(\Qstar)=\mathcal{N}({Q})$. 
\end{lemma}

\begin{example}\label{example coprime}\mbox{}
\begin{itemize} 
    \item[(a)] For MCAR$(p)$ models, $\Qstar(z)=I_k$ holds. Thus, $P(z)=\Pstar(z)$ and $Q(z)=I_k$ always provide a coprime right polynomial fractional description of $ H(z)$.
    \item[(b)]Consider an MCARMA$(2,1)$ model with coprime AR polynomial and MA polynomial given by
\begin{align*}
    \Pstar(z) = 
    \begin{pmatrix}
        (z+2)^2 & 0 \\
        0 & (z+2)^2
    \end{pmatrix} \quad \text{ and } \quad
    \Qstar(z) = 
    \begin{pmatrix}
        z+1 & 0 \\
        0 & z+1
    \end{pmatrix}.
\end{align*}    
Since $\Pstar(z)$ and $\Qstar(z)$ are diagonal matrix polynomials and are right coprime, the unique coprime right polynomial fraction description of the transfer function $ H(z)$ is given through $P(z)=\Pstar(z)$ and $Q(z)=\Qstar(z)$.

\item[(c)] Consider an MCARMA$(3,1)$ model with coprime AR polynomial and MA polynomial given by
\begin{align*}
    \Pstar(z) = 
    \begin{pmatrix}
        \frac{1}{4}(2z+3)(2z^2+7z+7) & -\frac{1}{4}(z+2)(3z+5) \\
        -(z+1)^2 & (z+1)^2(z+2)
    \end{pmatrix}  \, \text{ and } \,
    \Qstar(z) = - 
    \begin{pmatrix}
        z+1 & \frac{1}{4} \\
        0 & z+3
    \end{pmatrix}.
\end{align*}
Then the coprime right polynomial fractional description
of $ H(z)$ is given by
\begin{align*}
    P(z) = 
    \begin{pmatrix}
        (z+2)^3 & 0 \\
        0 & (z+1)^3
    \end{pmatrix} \quad \text{ and } \quad
    Q(z) = -
    \begin{pmatrix}
        z+2 & 1 \\
        1 & z+2
    \end{pmatrix}.
\end{align*}

    \item[(d)] The controller canonical state space model $(\BA, \BB, \BFC, L)$ in \Cref{example} with coprime right polynomial fraction description $P(z)$ and $Q(z)$ is as well an MCARMA$(2,1)$ model with coprime AR  and MA polynomial given by
    \begin{align*}
    \Pstar(z) = 
    \begin{pmatrix}
        z^2+z+1 & 0 & 0 \\
        0 & z^2+z+1 & 0 \\
        -\frac{1}{2} & \frac{1}{2} & z^2+z+1
    \end{pmatrix} \quad\text{ and } \quad
    \Qstar(z) = 
    \begin{pmatrix}
        z+1 & -z-1 & 0 \\
        z+1 & z+1 & 0 \\
        z+1 & z+1 & z+1
    \end{pmatrix}.
\end{align*}  
\end{itemize}
\end{example}

In particular, these examples further emphasise the existence of coprime right polynomial fraction descriptions  \eqref{eq: coprime right polynomial fraction description}  with polynomials as in \eqref{def PQ}. Since the aim of this paper is not the investigation of polynomial fraction descriptions, but the application of (local) orthogonality graphs to state space models, we do not investigate this further and move on to the topic of the invertibility of a state model.

\subsection{Invertible controller canonical state space models}\label{subsec: Inversion of MCARMA processes}
Suppose $\CY$ is a solution to a state space model that has a controller canonical representation $(\BA, \BB, \BFC, L)$, and we assume that the driving Lévy process satisfies the following common assumption.

\begin{assumptionletter}\label{Assumption on Levy process}
The $k$-dimensional  Lévy process $L=(L(t))_{t\in\R}$ satisfies \linebreak $\BE L(1)=0_k\in \R^k$ and  $\BE \Vert L(1) \Vert^2 < \infty$ with $\Sigma_L=\BE[L(1)L(1)^\top ]$.
\end{assumptionletter}
Then the second moments of $X(t)$ and thus of $Y(t)$ also exist \citep[Lemma A.4]{BR12}, which is a basic requirement for the forthcoming considerations on the existence of \mbox{(local) orthogonality graphs.}

Due to the state equation $Y(t) = \BFC X(t)$, we obtain the output process $Y(t)$ directly from the input process $X(t)$. However, the recovery of $X(t)$ from the output process $(Y(s))_{s\leq t}$, is not as obvious, because $\BFC$ is not invertible (cf. \Cref{Example 2.10}). Only for $q=0$, corresponding to an MCAR$(p)$ model, the simple structure of $\BFC$ allows to use the relation 
\begin{align} \label{block}
\begin{split}
    D^{(j-1)} Y(t)&= X^{(j)}(t), \quad\quad j=1,\ldots,p, \quad \quad \text{ where } \\
X^{(j)}(t)&=
\begin{pmatrix}
X_{(j-1)k+1}(t), & \cdots, & X_{jk}(t)
\end{pmatrix}^\top,
\end{split}
\end{align}
 is the \textsl{$j$-th $k$-block} of $X(t)$ and $D^{(1)} Y(t), \ldots, D^{(p-1)} Y(t)$ denote the mean-square derivatives of $Y(t)$. Therefore, in this case it is possible to recover  $X(t)$ from $Y(t)$ and its derivatives $D^{(1)} Y(t)$, \ldots, $D^{(p-1)} Y(t)$ via \eqref{block}. However, for controller canonical state space models with $q>0$, we cannot apply this approach. Indeed,  the structure of $\BA$ still yields
\begin{align} 
    D^{(1)}X^{(j)}(t)&=X^{(j+1)}(t), \hspace*{1.7cm} j=1,\ldots,p-1, \nonumber
\intertext{and together with $Y(t)=\BFC X(t) = C_0 X^{(1)}(t)+ \dots + C_q X^{(q+1)}(t)$, we obtain  that}
 \label{2.7}
    D^{(j-1)}Y(t)&= \sum_{i=0}^{q} C_{i} X^{(j+i)}(t), \quad \quad j=1,\ldots, p-q.
\end{align}
Consequently, the $p$ $k$-blocks of $X$ cannot generally be recovered from these $(p-q)$ equations. 

\pagebreak

\begin{remark}\label{Remark Derivatives in linear space}\mbox{}
\begin{itemize}
    \item[(a)] 
    For the reader's convenience, we define
    \begin{align}\label{eq: def C}
        \UBC \coloneqq (0_k, \ldots, 0_k, C_0, \ldots, C_q)\in M_{k\times kp}(\R).
    \end{align}
    From \eqref{2.7} we then receive the shorthands
    \begin{align} \label{eq: def CC}
        Y(t)=\OBC X(t) \quad \text{ and } \quad D^{(p-q-1)}Y(t)= \UBC X(t).
    \end{align}
     In particular, this implies that $\CY$ and its components $\CY_v$, $v\in V$, are $(p-q-1)$ times mean-square differentiable. 
      \item[(b)] A conclusion from (a) and  \cite{VF23pre}, Remark 2.6, is then that for $v\in V$ and $t\geq 0$, 
    \begin{align*}
        D^{(1)}Y_v(t), \dots, D^{(p-q-1)}Y_v(t) \in \mathcal{L}_{Y_v}(t).
    \end{align*}  
\end{itemize}
\end{remark}

For controller canonical state space models as in \eqref{CCform} with $q>0$, we overcome the challenge of recovering the state process from the output process under some mild assumptions. Of course, due to $q>0$, the class of MCAR$(p)$ models are excluded in the following considerations. However, this is not an essential limitation, because the (local) orthogonality graphs and the orthogonal projections, respectively for this case are already known \citep{VF23pre}. We first define causal \textsl{invertible} controller canonical state space models, which are a special subclass of controller canonical state space models. 


\begin{definition}\label{def: causal invertible state space model}
Let $(\Astar, \Bstar, \Cstar, L)$ be a state space model with controller canonical form $(\BA, \BB, \OBC, L)$ as in \eqref{CCform} and right coprime polynomials $P(z)$ and $Q(z)$ as in \eqref{def PQ} with $p>q>0$.
Suppose that
\begin{align}\label{Assumption Q}
    \rang(C_q)=k, \quad \mathcal{N}(Q)\subseteq (-\infty, 0) + i\R,  \quad \text{and} \quad \mathcal{N}(P)\subseteq (-\infty, 0) + i\R.
\end{align}
Then $(\BA, \BB, \OBC, L)$ is called a \textsl{causal invertible controller canonical state space (ICCSS) model} of order $(p,q)$ and the stationary solution $\CY=(Y(t))_{t\in\R}$ of the ICCSS$(p,q)$ model is called 
ICCSS$(p,q)$ process.
\end{definition}

\begin{remark}\mbox{}
\begin{itemize}
    \item[(a)]
 Since $\mathcal{N}({P}) = \sigma(\BA)\subseteq (-\infty, 0) + i\R$ \cite[Corollary 3.8]{MA07}, there exists a unique stationary solution $X(t)$ of the observation equation \cite[Theorem 4.1]{Sato:Yamazato} which has the representation
    \begin{align*} 
        X(t)=\int_{-\infty}^t e^{\BA (t-u)} \BB dL(u), \quad t\in \R.
    \end{align*}
Hence, there exists as well a stationary version of the output process $\CY$, which has the moving average representation
    \begin{align*}
        Y(t)=\uint g(t-u) dL(u) \quad \text{ with } \quad g(t)=  \OBC e^{\BA t} \BB {\boldsymbol 1}_{\{t\geq 0\}}, \quad t\in \R.
    \end{align*}
Throughout this paper, we are working with these stationary versions of $\CX$ and $\CY$. 
    \item[(b)] The assumptions on $Q(z)$ are necessary to recover $\CX$ from $\CY$ and to motivate the name ICCSS model, as we see in the remainder of this section.
    \item[(c)]  In the running \Cref{example}, we have $p=2>q=1>0$, $\rang(C_1)=3$, as well as $\mathcal{N}(Q) = \{-1\}  \subseteq (-\infty, 0) + i\R$, and $\mathcal{N}(P) = \sigma(\BA) = \{-\frac{1}{2} + \frac{\sqrt{3}}{2}i, -\frac{1}{2}- \frac{\sqrt{3}}{2}i  \}\subseteq (-\infty, 0) + i\R$. Thus, all of the assumptions in \eqref{Assumption Q} are satisfied. 
    $(\BA, \BB, \OBC, L)$ is a causal invertible controller canonical state space model of order $(2,1)$, and the stationary solution $\CY$ is an ICCSS$(2,1)$ process. Furthermore, the assumptions in \eqref{Assumption Q} are also satisfied in \Cref{example coprime}(b,c). 
\end{itemize}
\end{remark}

Under Assumption \eqref{Assumption Q}, \cite{BR12}, Lemma 4.1, derive a stochastic differential equation for the first $(kq)$ components of $\CX$. This follows simply from combining the first $q$ $k$-blocks of the state transition equation \linebreak $dX(t)= \BA X(t)dt+ \BB dL(t)$ with the observation equation $Y(t)=\OBC X(t)$ having the special structure of $\BA$, $\BB$, and $\OBC$ in mind. 

\begin{lemma}
Let $\CY$ be an 
ICCSS$(p,q)$ process with $p>q>0$.
Denote the $(kq)$-dimensional upper truncated state vector $\CX^q=(X^q(t))_{t\in \R}$ of $\CX$ by
\begin{align*}
X^q(t)=
\begin{pmatrix}
X^{(1)}(t)\\
\vdots\\
X^{(q)}(t)
\end{pmatrix}, \quad t\in \R,
\end{align*}
where $X^{(1)}(t),\ldots,X^{(q)}(t)$ are the $k$-dimensional random vectors as defined in \eqref{block}.
Then $\CX^q$ satisfies 
\begin{align}\label{dgl Xq}
dX^q(t)=\BLAM X^q(t)dt+ \BT Y(t)dt,
\end{align}
where $\sigma(\BLAM) \subseteq (-\infty, 0) + i\R$,
\begin{align*}
\BLAM =
\begin{pmatrix}
  0_k & I_k & 0_k & \cdots & 0_k \\
  0_k & 0_k & I_k & \ddots & \vdots \\
  \vdots &  & \ddots & \ddots & 0_k \\
  0_k & \cdots & \cdots & 0_k & I_k \\
  -C_q^{-1}C_0 & -C_q^{-1}C_1 & \cdots & \cdots & -C_q^{-1}C_{q-1}
\end{pmatrix}
\in M_{kq}(\R) \,\text{ and } \,
\BT=
\begin{pmatrix}
0_k \\
\vdots \\
0_k \\
C_q^{-1}
\end{pmatrix}
\in M_{kq\times k}(\R).
\end{align*}
\end{lemma}

\begin{remark}
\begin{itemize} $\mbox{}$
    \item[(a)] Assumption \eqref{Assumption Q} corresponds to the minimum-phase assumption in classical time series analysis \citep{HA12} and implies Assumption A2 in \cite{BR12}, who even allow for rectangular matrices $C_0$,\ldots,$C_q$. To see this, note that Assumption \eqref{Assumption Q} yields
    \begin{align*}
        \mathcal{N}(C_q^{-1} Q) &=  \{ z\in \C : \det(C_q^{-1} Q(z))=0\} =  \{ z\in \C : \det(Q(z))=0\} \\
        & = \mathcal{N}(Q) \subseteq (-\infty,0)+i\R,
    \end{align*}
    which is one of their assumptions. Furthermore, $\sigma(\BLAM)= \mathcal{N}(C_q^{-1} Q)$  (cf.~Lemma 3.8 by \citealp{MA07}). 
    Thus, $\BLAM$ has full rank and, due to the structure of $\BLAM$, we obtain that $C_q^{-1}C_0$ has full rank. It follows that $C_0$ and $(C_q)^\top C_0$ have full rank as well, which is the second assumption in \cite{BR12}.
    \item[(b)] If the AR polynomial $\Pstar(z)$ and the MA polynomial $\Qstar(z)$ 
    of an MCARMA model are left coprime,  Assumption \eqref{Assumption Q} can equally be made for  $\Pstar(z)$ and $\Qstar(z)$, respectively. Indeed, $\mathcal{N}(\Qstar)=\mathcal{N}(Q)$ and $\mathcal{N}(\Pstar) = \mathcal{N}(P)$ by \Cref{Lemma: Zusammenhang MCARMA}. Further, straightforward calculations of $\Qstar(z)P(z)=\Pstar(z)Q(z)$ give $\Qstar_q A_p = \Pstar_p C_q $. For $\mathcal{N}(\Pstar) = \mathcal{N}(P)\subseteq (-\infty,0)+i\R$ we have $0\notin\mathcal{N}(\Pstar)$ and thus, $\det (\Pstar_p)=\det (\Pstar(0))\not=0$. Similarly $\det (A_p) \not=0$ follows, so $\Pstar_p$ and $A_p$ are invertible. Hence, if $\mathcal{N}(\Pstar) = \mathcal{N}(P)\subseteq (-\infty,0)+i\R$, then $\Qstar_q$ has full rank, if and only if, $C_q$ has full rank.
\end{itemize}
\end{remark}

The differential equation \eqref{dgl Xq} has the solution
\begin{align}\label{erste komponenten}
X^q(t)=e^{\BLAM (t-s)}X^q(s)+\int_s^te^{\BLAM (t-u)} \BT Y(u)du, \quad  \: s<t,
\end{align}
\cite[(4.3)]{BR12}. Therefore, we can compute 
$X^q(t)$ based on the knowledge of the initial value $X^q(s)$ and $(Y( u ))_{s \leq u \leq t}$. In Propositions \ref{Xq als integral} and \ref{Xq als L2 integral} we even show the integral representation
\begin{align*}
X^q(t) = \int_{-\infty}^t e^{\BLAM (t-u)} \BT Y(u) du
\end{align*}
$\mathbb{P}$-a.s.~and in the mean square, respectively. Hence, $X^q(t)$ is even uniquely determined by the entire past $(Y(s))_{s\leq t}$, implying that the truncated state vector $\CX^q$ can be recovered from $\CY$. The remaining $k$-blocks $\CX^{(q+j)}$, $j=1,\ldots,p-q$, are obtained from $\CX^q$ and $\CY$ by differentiation as in \cite{BR12}, Lemma 4.2:

\begin{lemma}\label{restliche komponenten}
Let $\CY$ be an 
ICCSS$(p,q)$ process with $p>q>0$.
Then
\begin{align*}
X^{(q+j)}(t)
= \UBE^\top \Bigg[ \BLAM^j X^q(t) + \sum_{m=0}^{j-1}\BLAM^{j-1-m} \BT D^{(m)} Y(t)\Bigg], \quad j=1,\ldots,p-q, \: t\in \R.
\end{align*}
\end{lemma}

Note that there is a duplication of notation in \cite{BR12}, which can be seen by recalculating the induction start. We therefore give the corrected result in \Cref{restliche komponenten}.

\begin{example}
Coming back to \Cref{example}, the output process $Y$ has the representation as the linear combinations of $X$ via
     \begin{align*}
         Y(t) =  \OBC X(t)
            = \begin{pmatrix}
                X_1(t) - X_2(t) + X_4(t) - X_5(t) \\
                X_1(t) + X_2(t) + X_4(t) + X_5(t) \\
                X_1(t) + X_2(t) + X_3(t) + X_4(t) + X_5(t) +X_6(t)  
            \end{pmatrix}.
     \end{align*}
From this representation, it is not immediately obvious how $X$ can be recovered from $Y$. However, we can define
    \setlength{\arraycolsep}{4pt}
    \begin{align*}
    \BLAM = - C_1^{-1} C_0 =
        \begin{pmatrix*}[r]
        -1 & 0 & 0 \\
        0 & -1 & 0 \\
        0 & 0 & -1
        \end{pmatrix*} \quad \text{ and } \quad
    \BT = C_1^{-1} = 
        \frac{1}{2}\begin{pmatrix*}[r]
        1 & 1 & 0 \\
        - 1 & 1 & 0 \\
        0 & -2 & 2
        \end{pmatrix*},
    \end{align*}  
with $\sigma(\BLAM) =\mathcal{N}(Q) = \{-1\} \subseteq (-\infty, 0) + i\R$. Due to $\UBE^\top=I_3$ and the simple form of $\BLAM$ and its matrix exponential, respectively, we can recover the input process $X$ from the output process $Y$ by 
\begin{align*}
    X^{(1)}(t)
    &= X^q(t) 
    = \int_{-\infty}^t e^{-(t-u)} \frac{1}{2} 
    \begin{pmatrix} 
        \phantom{-} Y_1(u) + Y_2(u)  \\
        -Y_1(u) + Y_2(u) \\
        -2Y_2(u) + 2Y_3(u)
    \end{pmatrix} du,\\
    X^{(2)}(t) &= 
      - X^{(1)}(t) + \frac{1}{2} \begin{pmatrix} 
        \phantom{-} Y_1(t) + Y_2(t)  \\
        -Y_1(t) + Y_2(t) \\
        -2Y_2(t) + 2Y_3(t)
    \end{pmatrix}.
\end{align*}
\end{example}

In summary, in the example as well as in the general setting, we are able to compute not only the truncated state vector $X^q(t)$ but also the full state vector $X(t)$ based on the knowledge of $(Y(s))_{s \leq t}$. This justifies calling the ICCSS process $\CY$ \textsl{invertible} if Assumption \eqref{Assumption Q} holds.

\section{Orthogonal projections of ICCSS processes}\label{sec: Linear prediction of MCARMA processes}
In this section, we derive the orthogonal projections of ICCSS processes and their derivatives which we require to characterise (local) 
Granger causality and (local) contemporaneous correlation for ICCSS processes. First, we give alternative representations of $Y_a(t+h)$ as well as $D^{(p-q-1)}Y_a(t+h)$, $a\in V=\{1,\ldots,k\}$, suitable for the calculation of orthogonal projections in \Cref{sec: Representations of ICCSS processes and their derivatives}. Note that we consider the process $D^{(p-q-1)}Y_a(t+h)$ since, by \Cref{Remark highes derivative} below, it is the highest existing derivative of the ICCSS process which we require for the definition of local Granger causality and local contemporaneous correlation, respectively. In \Cref{sec: Orthogonal projections of ICCSS processes and their derivatives}, we then present the corresponding orthogonal projections of both random variables on $\mathcal{L}_{Y_S}(t)$ for $S\subseteq V$.   
Furthermore, we discuss the limit of the projections of difference quotients.

\subsection{Representations of ICCSS processes and their derivatives}\label{sec: Representations of ICCSS processes and their derivatives}
The aim of this subsection is to develop a $\mathbb{P}$-a.s.~representation of $Y_a(t+h)$ and $D^{(p-q-1)}Y_a(t+h)$, $a\in V$. Therefore, we first introduce the $\mathbb{P}$-a.s.~integral representation of the upper $q$-block truncation $\CX^q$, which is a multivariate generalisation of \cite{BR15}, Theorem 2.2. 

\begin{proposition}\label{Xq als integral}
Let $\CY$ be an 
ICCSS$(p,q)$ process with $p>q>0$. 
 Then, for all $t\in \R$, we have
\begin{align*}
  X^q(t) = \int_{-\infty}^t e^{\BLAM (t-u)} \BT Y(u) du \quad \text{$\mathbb{P}$-a.s.}
\end{align*}
\end{proposition}

Due to the well-definedness of this integral, it is obvious that the following representations of $\CY$ and its derivatives are well-defined as well. 

\begin{theorem}\label{Darstelung von Yb}
Let $\CY$ be an 
ICCSS$(p,q)$ process with $p>q>0$.
Then, for $h\geq 0$, $t\in \R$, and $a\in V$, it holds that
\begin{align*}
Y_a(t+h)
= & \: \int_{-\infty}^t  e_a^\top \OBM(h) e^{\BLAM (t-u)} \BT Y(u) du \\
&+ \sum_{m=0}^{p-q-1} e_a^\top \OBMnu(h) \BT D^{(m)} Y(t) +  e_a^\top \Oeps(t,h) \quad \text{$\mathbb{P}$-a.s.}
\quad  \text{and} \\
D^{(p-q-1)}Y_a(t+h) 
= & \: \int_{-\infty}^t  e_a^\top \BM(h) e^{\BLAM (t-u)} \BT Y(u) du \\
&+ \sum_{m=0}^{p-q-1} e_a^\top \BMnu(h) \BT D^{(m)} Y(t) +  e_a^\top \Ueps(t,h) \quad \text{$\mathbb{P}$-a.s.}
\end{align*}
Here, we abbreviate 
\begin{gather*}
\OBM(h) = \OBC e^{\BA h} \Bigg( \OBE + \sum_{j=1}^{p-q} \BFE_{q+j} \UBE^\top \BLAM^j \Bigg),\quad\quad
\BM(h) = \UBC e^{\BA h} \Bigg( \OBE + \sum_{j=1}^{p-q} \BFE_{q+j} \UBE^\top \BLAM^j \Bigg),\\
\OBMnu(h) = \OBC e^{\BA h} \sum_{j=m + 1}^{p-q} \BFE_{q+j} \UBE^\top \BLAM^{j-1-m},\quad\quad
\BMnu(h) = \UBC e^{\BA h} \sum_{j=m + 1}^{p-q} \BFE_{q+j} \UBE^\top \BLAM^{j-1-m},\\
 \Oeps(t,h) =  \OBC \int_t^{t+h} e^{\BA (t+h-u)} \BB dL(u), \quad\quad
 \Ueps(t,h) =  \UBC \int_t^{t+h} e^{\BA (t+h-u)} \BB dL(u).
\end{gather*}
where $\UBC$ is defined in \eqref{eq: def C}, $\OBE$ and $\UBE$ are defined in \eqref{eq: def E}, and $\BFE_j$ is defined in \eqref{eq: def Ej}, $j=1,\ldots,p$. Finally, $\Oeps(t,0)=\Ueps(t,0)=0_k \in \R^k$.
\end{theorem}


\begin{remark}\label{remark: comparison to MCAR 1}
For an MCAR$(p)$ process, \cite{VF23pre} state in Lemma 6.8 that
\begin{align} \label{eqaa}
Y_a(t+h)&= e_a^\top \OBC e^{\BA h}  \sum_{m=0}^{p-1} \BFE_{m+1} D^{(m)} Y(t) + e_a^\top \OBC \int_t^{t+h} e^{\BA (t+h-u)} \BB dL(u) \quad \text{$\mathbb{P}$-a.s.}
\end{align}
Thus, if we want to compare our \Cref{Darstelung von Yb} to the results for MCAR$(p)$ processes, we have to interpret 
\begin{align*}
    \OBMnu(h) \BT \: \widehat{=} \: \OBC e^{\BA h}\BFE_{m+1}, \quad m=0,\ldots,p-1, \quad \text{ and } \quad 
    \OBM(h) e^{\BLAM (t-u)} \BT \: \widehat{=} \: 0_{k} \text{ for } u<t.
\end{align*}
Then \eqref{eqaa} can be seen as a special case of \Cref{Darstelung von Yb}. Let us briefly heuristically justify that this interpretation is reasonable. First of all, in $\OBMnu(h)\BT$ the summand $j=m+1$ is mainly relevant. For this summand we have with $\BLAM^0=I_{kq}$ that
\begin{align}\label{eq: Rest a}
    \OBC e^{\BA h} \BFE_{q+m+1} \UBE^\top \BT
    = \OBC e^{\BA h} \BFE_{q+m+1} C_q^{-1}.
\end{align}
If $q=0$ is inserted into $\OBMnu(h)\BT$, all summands disappear due to the zero dimensionality of $\BLAM^{j-1-m}$, $j=m+2,\ldots,p-q$, except for \eqref{eq: Rest a}. With $C_q=I_k$ it remains as claimed $ \OBMnu(h) \BT  	\: \widehat{=} \: \OBC e^{\BA h} \BFE_{m+1}$ for $m=0,\ldots,p-1$. For the second matrix function $\OBM(h) e^{\BLAM (t-u)} \BT$, $u<t$, we use similar arguments to show that it can be interpreted as a zero matrix. Although we get a non-zero matrix for $t=u$, this event is a Lebesgue null-set.

\end{remark}

\begin{example}\label{Example Darstellung}
In \Cref{example}, we have $m=p-q-1=0$, so 
\begin{align*}
Y(t+h) = D^{(0)}Y(t+h)
= \int_{-\infty}^t  \OBM(h) e^{\BLAM (t-u)} \BT Y(u) du 
+ \OBMo(h) \BT Y(t) +  \Oeps(t,h) \quad \text{$\mathbb{P}$-a.s.}
\end{align*}
If we abbreviate $c(h) := 3 \cos(\sqrt{3} h/ 2)$ and $s(h): = \sqrt{3}\sin (\sqrt{3} h/ 2)$, the three addends can be specified as follows. 
\begin{align*}
    & \int_{-\infty}^t  \OBM(h) e^{\BLAM (t-u)} \BT Y(u) du \\
    & \quad = \frac{e^{-\frac{h}{2}}}{3} \int_{-\infty}^t e^{-(t - u)} 
    \left(\begin{array}{c}
        -2 s(h) Y_1(u)
            \phantom{\: \: + \frac{1}{3}\left[h c(h)+s(h)\right]\left( Y_1(u) - Y_2(u) \right)}\\
        -2 s(h) Y_2(u)
            \phantom{\: \: + \frac{1}{3}\left[h c(h)+s(h)\right]\left( Y_1(u) - Y_2(u) \right)}\\
        -2 s(h) Y_3(u)
            + \frac{1}{3}\left[h c(h)+s(h)\right]\left( Y_1(u) - Y_2(u) \right)  
\end{array}\right)  du, \\
   & \OBMo(h) \BT Y(t) \\
   & \quad =  \frac{e^{-\frac{h}{2}}}{3} \begin{pmatrix}
    \left[ c(h) + s(h) \right] Y_1(t) 
        \phantom{\:\: +\frac{1}{6}\left[- h c(h) + \left(3h+2\right)s(h) \right] \left(Y_1(t) - Y_2(t) \right)}\\
    \left[ c(h) + s(h) \right] Y_2(t)
        \phantom{\:\: +\frac{1}{6}\left[- h c(h) + \left(3h+2\right)s(h) \right] \left(Y_1(t) - Y_2(t) \right)}\\
    \left[ c(h) + s(h) \right] Y_3(t)
        +\frac{1}{6}\left[- h c(h) + \left(3h+2\right)s(h) \right] \left(Y_1(t) - Y_2(t) \right)
    \end{pmatrix},\\
  & \Oeps(t,h) \\
 & \quad = \frac{e^{-\frac{h}{2}}}{3} \int_{t}^{t+h} e^{-\frac{(t-u)}{2}} 
  \left(\begin{matrix*}[c]
    \left[ c(t + h - u) + s(t + h - u) \right] \left( dL_1(u) - dL_2(u) \right) \\
    \left[ c(t + h - u) + s(t + h - u) \right] \left( dL_1(u) + dL_2(u) \right) \\
    \left[ c(t + h - u) + s(t + h - u) \right] \left( dL_1(u) + dL_3(u) \right)
  \end{matrix*}\right)\\
  & \quad\quad + 
  \begin{pmatrix}
  \: 0\\
   \: 0\\
\frac{e^{-\frac{h}{2}}}{3} \int_{t}^{t+h} e^{-\frac{(t-u)}{2}} \left[\left( 1+\frac{t+h-u}{3} \right) c(t + h - u)
+ \left( \frac{1}{3} - \left(t+h-u \right) \right) s(t+h-u)\right] dL_2(u)
  \end{pmatrix} .
\end{align*}

\end{example}

\subsection{Orthogonal projections of ICCSS processes and their derivatives}\label{sec: Orthogonal projections of ICCSS processes and their derivatives}
The representations of $Y_a(t+h)$ and $D^{(p-q-1)}Y_a(t+h)$ in \Cref{Darstelung von Yb} suggest that for the orthogonal projection of the $a$-th component one time step 
into the future, on the one hand, the past $(Y_V(s))_{s\leq t}$ of all components and on the other hand, the future of the Lévy process $(L(s)-L(t))_{t \leq s \leq t+h}$ is relevant. However, for a formal proof, we require that all integrals are defined in $L^2$. Therefore, we show that the integral representation of $\CX^q$ in \Cref{Xq als integral} holds in $L^2$. The proof is based on the ideas of the proof of  Theorem 2.8 in \cite{BR15}  in the univariate setting.

\begin{proposition} \label{Xq als L2 integral}
Let $\CY$ be an 
ICCSS$(p,q)$ process with $p>q>0$.
Then, for \mbox{$a,v\in V$} and $t\in \R$, the integral
\begin{align*}
\int_{-\infty}^t e_a^\top e^{\BLAM (t-u)} \BT e_v Y_v(u) du\in \mathcal{L}_{Y_v}(t)
\end{align*}
exists as $L^2$-limit. In particular,
$
X^q(t) = \int_{-\infty}^t e^{\BLAM (t-u)} \BT Y(u) du
$
exists as $L^2$-limit. 
\end{proposition}


Before finally moving on to the orthogonal projections, we introduce one last alternative representation, this time for the difference quotient $(D^{(p-q-1)} Y_a(t+h)- D^{(p-q-1)} Y_a(t))/h$.  With this representation we can argue that $D^{(p-q-1)} Y_a(t)$ is indeed the maximum derivative of $Y_a(t)$ which we need for local Granger causality and local contemporaneous correlation.

\begin{lemma} \label{Lemma 4.7}
Let $\CY$ be an 
ICCSS$(p,q)$ process with $p>q>0$.
Then for $h \geq 0$, $t \in \R$, and $a\in V$ the representation
\begin{align*}
&\frac{D^{(p-q-1)} Y_a(t+h)- D^{(p-q-1)} Y_a(t)}{h} \\
&\quad \quad =   \int_{-\infty}^t  e_a^\top \BM'(0) e^{\BLAM (t-u)} \BT Y(u) du 
+ \sum_{m=0}^{p-q-1} e_a^\top \BMnu'(0) \BT D^{(m)} Y(t) \\ 
& \quad\quad \quad+ e_a^\top O(h)\bR_1+ e_a^\top O(h)\bR_2  + e_a^\top\frac{\Ueps(t,h)}{h}
\end{align*}
holds, where $\bR_1,\bR_2$ are random vectors in $\mathcal{L}_{Y}(t)\subseteq L^2$. $\BM'(0)$ and $\BMnu'(0)$ denote the first derivatives of $\BM(\cdot)$ and $\BMnu(\cdot)$ in zero. The random variable $e_a^\top \Ueps(t,h)/h$
is independent of the former summands and 
\begin{eqnarray*}
    \lim_{h \downarrow 0}\frac{1}{h}\BE\left[(e_a^\top\Ueps(t,h))^2\right]= e_a^\top \UBC \BB \Sigma_L \BB^\top \UBC e_a\not=0
    \quad \text{ but }  \quad \lim_{h\downarrow 0}\frac{1}{h^2}\BE\left[(e_a^\top\Ueps(t,h))^2\right]=\infty.
\end{eqnarray*}
\end{lemma}

\begin{remark}\label{Remark highes derivative}
An important consequence of \Cref{Lemma 4.7} is that the mean-square limit of the difference quotient does not exist, and hence, for all components of the ICCSS process no mean-square derivatives higher than $(p-q-1)$ exist. Thus, for local Granger causality and local contemporaneous correlation, we must always analyse the $(p-q-1)$-th derivative. It also becomes clear that in the definition of local contemporaneous correlation, one must divide by $h$ and not by $h^2$.
\end{remark}


Now, we specify the orthogonal projections. 

\begin{theorem}\label{Projektionen für MCARMA}
Let $\CY$ be an 
ICCSS$(p,q)$ process with $p>q>0$. 
Suppose $S \subseteq V$ and $a\in V$. Then, for $h\geq 0$ and $t\in \R$, we have
\begin{align*}
P_{\mathcal{L}_{Y_S}(t)} Y_a(t+h)
= &\sum_{v \in S} \int_{-\infty}^t  e_a^\top \OBM(h)  e^{\BLAM (t-u)} \BT e_v Y_v(u) du \\
  & + \sum_{v \in S} \sum_{m=0}^{p-q-1} e_a^\top\OBMnu (h) \BT e_v D^{(m)} Y_v(t) \\
  & + P_{\mathcal{L}_{Y_S}(t)} \Bigg( \sum_{v \in V\setminus S} \int_{-\infty}^t e_a^\top\OBM(h) e^{\BLAM (t-u)} \BT e_v Y_v(u) du \Bigg) \\
  & + P_{\mathcal{L}_{Y_S}(t)} \Bigg(\sum_{v \in V\setminus S} \sum_{m=0}^{p-q-1} e_a^\top\OBMnu (h) \BT e_v D^{(m)} Y_v(t) \Bigg) \quad \text{$\mathbb{P}$-a.s.}
  \intertext{and}
P_{\mathcal{L}_{Y_S}(t)} D^{(p-q-1)}Y_a(t+h)
= &\sum_{v \in S} \int_{-\infty}^t  e_a^\top \BM(h)  e^{\BLAM (t-u)} \BT e_v Y_v(u) du \\
  & + \sum_{v \in S} \sum_{m=0}^{p-q-1} e_a^\top\BMnu (h) \BT e_v D^{(m)} Y_v(t) \\
  & + P_{\mathcal{L}_{Y_S}(t)} \Bigg( \sum_{v\in V\setminus S} \int_{-\infty}^t e_a^\top\BM(h) e^{\BLAM (t-u)} \BT e_v Y_v(u) du \Bigg) \\
  & + P_{\mathcal{L}_{Y_S}(t)} \Bigg(\sum_{v \in V\setminus S} \sum_{m=0}^{p-q-1} e_a^\top\BMnu (h) \BT e_v D^{(m)} Y_v(t) \Bigg) \quad \text{$\mathbb{P}$-a.s.},
\end{align*}
where $\OBM(\cdot)$, $\BM(\cdot)$, $\OBMnu (\cdot)$ and $\BMnu(\cdot)$ are defined in \Cref{Darstelung von Yb}.
\end{theorem}

The basic idea of the proof is simple: In the representation in \Cref{Darstelung von Yb}, the terms  $Y_a(t)$, its derivatives and integrals over the past are already in the linear space $\mathcal{L}_{Y_S}(t)$ if $a\in S$ (\Cref{Remark Derivatives in linear space} and \Cref{Xq als L2 integral}) and are therefore projected onto themselves. Furthermore, $(Y_S(s))_{s\leq t}$  and $(L(s)-L(t))_{t \leq s \leq t+h}$ are independent such that $e_a^\top \Oeps(t,h)$ and $e_a^\top \Ueps(t,h)$, respectively, are independent of $\mathcal{L}_{Y_S}(t)$ and are projected onto zero. Of course, this argument can also be used in \Cref{example} (\Cref{Example Darstellung}) to display the desired projections directly, and we refrain from specifying them.

\begin{remark}
When calculating the orthogonal projections, it becomes clear why we require  Assumption \eqref{Assumption Q}, a sufficient assumption to recover $X(t)$ from $(Y(s))_{s\leq t}$. Only then are we able to project the input process $X(t)$ onto the linear space of the output process $\mathcal{L}_{Y_S}(t)$.
\end{remark}

To apply local Granger causality and local contemporaneous correlation to 
ICCSS processes, we also need the following orthogonal projections. 

\begin{theorem}\label{Projektionen mit limneu}
Let $\CY$ be an 
ICCSS$(p,q)$ process with $p>q>0$. 
Suppose $S \subseteq V$ and $a\in V$. Then, for $t\in \R$, we have
\begin{align*}
& \limh P_{\mathcal{L}_{Y_{S}}(t)} \left(\frac{D^{(p-q-1)} Y_a(t+h)- D^{(p-q-1)} Y_a(t)}{h}\right) \\
&\quad =\: 
\sum_{v\in S} \int_{-\infty}^t e_a^\top  \BM'(0) e^{\BLAM(t-u)} \BT e_v Y_v(u) du 
+ \sum_{v \in S} \sum_{m=0}^{p-q-1} e_a^\top \BMnu'(0) \BT e_v D^{(m)} Y_v(t) \\
&\quad\phantom{=} 
+ P_{\mathcal{L}_{Y_{S}}(t)} \Bigg(\sum_{v\in V\setminus S} \int_{-\infty}^t e_a^\top \BM'(0) e^{\BLAM(t-u)} \BT e_v Y_v(u) du \Bigg) \\
&\quad\phantom{=} 
+ P_{\mathcal{L}_{Y_{S}}(t)} \Bigg(\sum_{v \in V\setminus S} \sum_{m=0}^{p-q-1} e_a^\top \BMnu'(0) \BT e_v D^{(m)} Y_v(t) \Bigg) 
\quad \mathbb{P}\text{-a.s.}
\end{align*}
and for $h\geq 0$,
\begin{align*}
 D^{(p-q-1)} Y_a(t+h)- P_{\mathcal{L}_{Y}(t)} D^{(p-q-1)} Y_a(t+h)
 = e_a^\top \Ueps(t,h)
\quad \mathbb{P}\text{-a.s.},
\end{align*}
where $\BM(\cdot)$, $\BMnu(\cdot)$, and $\Ueps(\cdot,\cdot)$ are defined in \Cref{Darstelung von Yb}. 
\end{theorem}

In this paper, for the derivation of the (local) orthogonality graph for ICCSS processes, the special case $S=V$ is most relevant, where a few terms are simplified.

\begin{corollary}\label{Projektion für MCARMA für S=V} 
Let $\CY$ be an 
ICCSS$(p,q)$ process with $p>q>0$.
Then, for $t\in\R$,  $h\geq 0$, and $a\in V$ the following projections hold.
\begin{itemize}
        \item[(a)] ${\displaystyle P_{\mathcal{L}_{Y}(t)}Y_a(t+h)
= e_a^\top \OBC e^{\BA h} X(t)}$ \quad \text{$\mathbb{P}$-a.s.},
\item[(b)]
${\displaystyle P_{\mathcal{L}_{Y}(t)} D^{(p-q-1)}Y_a(t+h)
= e_a^\top \UBC e^{\BA h} X(t)}$ \quad \text{$\mathbb{P}$-a.s.},
\item[(c)]
${\displaystyle \limh P_{\mathcal{L}_{Y}(t)} \left(\frac{D^{(p-q-1)} Y_a(t+h)- D^{(p-q-1)} Y_a(t)}{h}\right)
= e_a^\top \UBC \BA X(t) \quad \text{$\mathbb{P}$-a.s.}}$
\end{itemize}
 \end{corollary}
 
From \Cref{Projektion für MCARMA für S=V}(c), not only the existence of the limit  becomes clear, but also that of the limit
     \begin{align*}
         \limh P_{\mathcal{L}_{Y_{S}}(t)} \left(\frac{D^{(p-q-1)} Y_a(t+h)- D^{(p-q-1)} Y_a(t)}{h}\right) 
         = P_{\mathcal{L}_{Y_{S}}(t)} \left( e_a^\top \UBC \BA X(t)\right).
     \end{align*}
The existence of these limits is essential for the well-definedness of local Granger causality and local contemporaneous correlation for ICCSS processes. 

\begin{remark} $\mbox{}$
\begin{itemize}
    \item[(a)] Although the derivation of the orthogonal projections for MCAR$(p)$ processes differs from that for ICCSS$(p,q)$ processes with $q>0$, the results are consistent with \cite{VF23pre}, Proposition 6.9 and Lemma 6.11 for MCAR processes, if we interpret $\OBM(h) e^{\BLAM (t-u)} \BT \: \widehat{=} \: 0_{k}$ for $u<t$ and $\OBMnu(h) \BT \:\widehat{=}\: \BFC e^{\BA h}\BFE_{m}$ as in \Cref{remark: comparison to MCAR 1}.
     \item[(b)] The linear projections in \Cref{Projektion für MCARMA für S=V}(a) match the linear projections for univariate CARMA processes in \cite{BR15}, Theorem 2.8. \cite{BA19} derives as well linear projections for MCARMA processes, but the results there differ from \cite{BR15}.
 \end{itemize}
\end{remark}

\section{Orthogonality graphs for ICCSS processes}\label{sec: orthogonality graph for MCARMA processes}

In this section, we derive (local) orthogonality graphs for state space models and obtain as the main result the characterisation of the directed and the undirected edges of the (local) orthogonality graph by the model parameters of the unique controller canonical form if this is an ICCSS$(p,q)$ model with $p>q>0$. 
To define the (local) orthogonality graph for ICCSS processes according to \Cref{Definition orthogonality graph}, certain requirements for the well-definedness must be met. We have already assumed that we use the stationary version of the ICCSS process throughout the paper and it has expectation zero. Furthermore, the continuity in the mean square of an ICCSS process is well known, it follows directly from \cite{CR39}, Lemma 1, since the covariance function is continuous in $0$. Therefore, we only need to make sure that the Assumptions \ref{Assumption an Dichte} and \ref{Assumption purely nondeterministic of full rank} are satisfied.

\begin{theorem} \label{Graph well defined}
Let $\CY$ be an $k$-dimensional  
ICCSS$(p,q)$ process with \mbox{$\Sigma_L>0$.} Then $\CY$ satisfies Assumptions \ref{Assumption an Dichte} and \ref{Assumption purely nondeterministic of full rank} and thus, the orthogonality graph and the local orthogonality graph are well defined and the Markov properties in \Cref{rem:Markov properties} hold.
\end{theorem}

\begin{remark}\mbox{}
\begin{itemize}
    \item[(a)] In principle, more general state space models $(\Astar, \Bstar, \Cstar, L)$  also satisfy  \Cref{Assumption graph}. The proof of \Cref{Graph well defined} shows that sufficient assumptions for the stationary state space process are that the driving Lévy process satisfies \Cref{Assumption on Levy process}, $\sigma(\Astar)\subseteq (-\infty,0)+i\R$, $f_{YY}(\lambda)>0$ $\forall \: \lambda\in\R$,  and $\Cstar \Bstar \BS_L \Bstar^{ \top} \Cstar^{\top} >0$. Then the (local) orthogonality graphs are well defined as well, and the Markov properties in \Cref{rem:Markov properties} hold. However, in this general context, we are not able to calculate the orthogonal projections needed to characterise the edges, which is our main goal. 
    \item[(b)] In our running \Cref{example}, we already know that $Y$ is an ICCSS(2,1) process and, since $\Sigma_L=I_3>0$, the orthogonality graph and the local orthogonality graph are well defined and the Markov properties hold.
\end{itemize}

\end{remark}

Let us now focus on the main results, i.e.,~the characterisations for the directed and the undirected edges in the (local) orthogonality graph for ICCSS processes. 
First, we present the characterisations of the (local) Granger non-causality.

\begin{theorem}\label{Zweite Charakterisierung gerichtete Kante für MCARMA}
Let $\CY=\CY_V$ be an 
ICCSS$(p,q)$ process with $p>q>0$ and $\Sigma_L>0$.
Let $a,b\in V$ and $a \neq b$. Then the following characterisations hold:
\begin{itemize}
    \item[(a)]{\makebox[2.7cm][l]{$\CY_a \nrarrow \CY_b \: \vert \: \CY_V$} 
    $\Leftrightarrow \quad
    e_b^\top \OBC \BA^\alpha \big( \OBE + \sum_{j=1}^{p-q} \BFE_{q+j} \UBE^\top \BLAM^j \big) \BLAM^\beta \BT e_a = 0$ \quad \text{and} \quad \\
    \phantom{a}\hspace{3.9cm}
    $e_b^\top \OBC \BA^\alpha \big( \sum_{j=m + 1}^{p-q} \BFE_{q+j} \UBE^\top \BLAM^{j-1-m}  \big) \BT e_a = 0,$ \\
    \phantom{a}\hspace{3.9cm}
    for $\alpha=0,\ldots,kp-1$, $\beta=0,\ldots,kq-1$, $m=0,\ldots,p-q-1$.}
    \item[(b)]{\makebox[2.7cm][l]{$ \CY_a \nrarrownull \CY_b \: \vert \: \CY_V $} 
    $\Leftrightarrow \quad
    e_b^\top \UBC \BA \big( \OBE + \sum_{j=1}^{p-q} \BFE_{q+j} \UBE^\top \BLAM^j \big) \BLAM^\beta \BT e_a = 0$ \quad \text{and} \quad \\
    \phantom{a}\hspace{3.9cm}
    $e_b^\top \UBC \BA \big( \sum_{j=m + 1}^{p-q} \BFE_{q+j} \UBE^\top \BLAM^{j-1-m}  \big) \BT e_a = 0,$\\
    \phantom{a}\hspace{3.9cm} for $\beta=0,\ldots,kq-1$, $m=0,\ldots,p-q-1$.}
    \end{itemize}
\end{theorem}

The basis for the proof of \Cref{Zweite Charakterisierung gerichtete Kante für MCARMA} are the following characterisations of the directed edges in  \Cref{Erste Charakterisierung gerichtete Kante für MCARMA}. The characterisations in \Cref{Erste Charakterisierung gerichtete Kante für MCARMA} are in turn developed from the definition of the directed edges in \Cref{subsec: orthogonality graphs} and the orthogonal projections of the ICCSS process and its derivatives in \Cref{sec: Linear prediction of MCARMA processes}.

\begin{proposition}\label{Erste Charakterisierung gerichtete Kante für MCARMA}
Let $\CY=\CY_V$ be an 
ICCSS$(p,q)$ process with $p>q>0$ and \mbox{$\Sigma_L>0$.}
Let $a,b\in V$ and $a \neq b$. Then the following characterisations hold:
\begin{itemize}
    \item[(a)]{\makebox[2.7cm][l]{$\CY_a \nrarrow \CY_b \: \vert \: \CY_V$} 
    $\Leftrightarrow \quad
    e_b^\top \OBM(h) e^{\BLAM t} \BT e_a = 0 \quad \text{and} \quad
    e_b^\top \OBMnu (h) \BT e_a = 0,$ \\
    \phantom{a}\hspace{3.9cm}
    for $m=0,\ldots,p-q-1$, $0 \leq h \leq 1$, $t \geq 0$.}
    \item[(b)]{\makebox[2.7cm][l]{$ \CY_a \nrarrownull \CY_b \: \vert \: \CY_V $} 
    $\Leftrightarrow \quad
    e_b^\top \BM'(0) e^{\BLAM t} \BT e_a = 0 \quad \text{and} \quad e_b^\top \BMnu'(0) \BT e_a = 0,$ \\
    \phantom{a}\hspace{3.9cm}
    for $m=0,\ldots,p-q-1$, $t \geq 0$.}
\end{itemize}
\end{proposition}

\begin{remark}\mbox{}
\begin{enumerate}
    \item[(a)] Except for the assumption $\Sigma_L>0$, the characterisations of the directed edges in the (local) orthogonality graph are independent of the chosen Lévy process, which is quite surprising. Thus, for example, the characterisations of (local) Granger causality  are the same for a Brownian motion driven ICCSS process and a Poisson driven ICCSS process with the same controller canonical triple $(\BA, \BB, \OBC)$, even if they have different covariance matrices and even though the path properties of these processes are significantly different. Note that in \Cref{example} we did not specify which Lévy process we are using.
    \item[(b)] The characterisation in \Cref{Erste Charakterisierung gerichtete Kante für MCARMA}(a) seems to depend on $h$. However, this is not the case as can be seen from \Cref{Zweite Charakterisierung gerichtete Kante für MCARMA}(a). So it does not matter whether we define directed edges by looking at the period $0\leq h \leq 1$ or by looking at the entire future $h\geq 0$. In terms of \cite{VF23pre}, there is no difference between Granger causality and global Granger causality for ICCSS processes.
\end{enumerate}
\end{remark}

Next, we present the characterisations of the undirected edges, i.e.,~contemporaneous uncorrelatedness.


\begin{proposition}\label{Einfache Charakterisierung ungerichtete Kante für MCARMA}
Let $\CY=\CY_V$ be an 
ICCSS$(p,q)$ process with $p>q>0$ and \mbox{$\Sigma_L>0$.} 
Let $a,b\in V$ and $a \neq b$. Then the following characterisations hold:
\begin{itemize}
    \item[(a)]\makebox[2.5cm][l]{$\CY_a \nsim \CY_b \: \vert \: \CY_V$} 
    $\Leftrightarrow \quad
    e_a^\top \int_0^{\min(h,\widetilde{h})} \OBC e^{\BA (h-s)} \BB \BS_L \BB^\top e^{\BA^\top(\widetilde{h}-s)} \OBC^\top ds \: e_b =0,$\\
    \phantom{a}\hspace{3.7cm}
    for $0\leq h, \widetilde{h}\leq 1$.\\
    \phantom{a} \hspace{2.65cm} 
    $\Leftrightarrow \quad e_a^\top \OBC \BA^\alpha \BB \BS_L \BB^\top \big(\BA^\top \big)^\beta \OBC^\top e_b=0,$ \\
    \phantom{a}\hspace{3.7cm}
    for $\alpha,\beta=0,\ldots,kp-1$.
    \item[(b)]{\makebox[2.5cm][l]{$ \CY_a \nsimnull \CY_b \: \vert \: \CY_V $} 
    $\Leftrightarrow \quad
    e_a^\top \UBC \BB \Sigma_L \BB^\top \UBC^\top e_b 
    = e_a^\top C_q \Sigma_L C_q^\top e_b =0.$}
\end{itemize}
\end{proposition}

The proof of this result again uses the orthogonal projections of  ICCSS processes and its derivatives of \Cref{sec: Linear prediction of MCARMA processes} and the definition of undirected edges of \Cref{subsec: orthogonality graphs}. The assumption $\Sigma_L>0$ is only used for the second characterisation in \Cref{Einfache Charakterisierung ungerichtete Kante für MCARMA}(a). 
However, it was also important for the proof of \Cref{Assumption graph}.

\begin{remark} \mbox{}
\begin{itemize}
    \item[(a)] The characterisations and thus, the undirected edges in the (local) orthogonality graph depend on the chosen Lévy process only by $\Sigma_L$. For example, the characterisations of the (local) contemporaneous correlation and thus the (local) orthogonality graph are the same for a Brownian motion driven ICCSS process and a Poisson driven ICCSS process with the same controller canonical triple $(\BA, \BB, \OBC)$ if both Lévy processes have the same covariance matrix $\Sigma_L$. However, in contrast to (local) Granger causality, it is necessary that the Brownian motion and the Poisson process have the same covariance matrix. Again, in our running \Cref{example}, we did not specify which Lévy process we were using, but we did specify that $\Sigma_L=I_3$.
    \item[(b)] The second characterisation in \Cref{Einfache Charakterisierung ungerichtete Kante für MCARMA}(a) shows that there is indeed no dependence on the lag $h$ again. As for the directed edges, it does not matter whether we define undirected edges by looking at the period $0\leq h, \widetilde{h} \leq 1$ or by looking at the entire future $h, \widetilde{h} \geq 0$. In terms of \cite{VF23pre}, there is no difference between contemporaneous correlation and global contemporaneous correlation for ICCSS processes. 
\end{itemize}
\end{remark}

We make some further comments on the characterisations in Propositions \ref{Erste Charakterisierung gerichtete Kante für MCARMA} and \ref{Einfache Charakterisierung ungerichtete Kante für MCARMA}. In particular, we compare the characterisations with each other and with the results in the literature, additionally we give some interpretations.

\begin{remark} $\mbox{}$
\begin{itemize}
    \item[(a)]
     The uniqueness of the polynomials $P(z)$ and $Q(z)$ in \eqref{eq: coprime right polynomial fraction description} (see \Cref{lemma: uniqueness}) leads to the uniqueness of the controller canonical state space representation, which in turn leads to the uniqueness of the edges in the (local) orthogonality graph. 
    \item[(b)] It can be shown by a simple calculation that $\OBC \BA^{p-q} = \UBC \BA $. Comparing \Cref{Zweite Charakterisierung gerichtete Kante für MCARMA}(a) and (b), we receive that Granger non-causality implies local Granger non-causality, which we know as well from the theory in \cite{VF23pre}. Similarly, $\OBC \BA^{p-q-1} = \UBC$, so comparing \Cref{Einfache Charakterisierung ungerichtete Kante für MCARMA}(a) and (b), we get that contemporaneous uncorrelatedness implies local contemporaneous uncorrelatedness, which is again in agreement with the theory. The relationships between Granger non-causality and local Granger non-causality, as well as contemporaneous uncorrelatedness and local contemporaneous uncorrelatedness in a general setting, are discussed in \cite{VF23pre}, Lemma 3.13 and Lemma 4.10.
    \item[(c)] In \Cref{example}, straightforward computations yield the (local) Granger causality relations and (local) contemporaneous correlations, visualised in the corresponding (local) causality graphs in \Cref{fig: Two graphical models b}. Again, the relationships between Granger non-causality and local Granger non-causality, and between contemporaneous uncorrelation and local contemporaneous uncorrelation are evident.
\end{itemize}
\end{remark}


\begin{interpretation}[Orthogonality graph]
To interpret the directed and the undirected edges in the orthogonality graph $G_{OG}$, we 
recall the representation of the $b$-th component 
\begin{align*}
Y_b(t+h)
= & \: \int_{-\infty}^t e_b^\top \OBM(h) e^{\BLAM (t-u)} \BT Y_V(u) du + \sum_{m=0}^{p-q-1} e_b^\top \OBMnu (h) \BT D^{(m)} Y_V(t) +  e_b^\top \Oeps(t,h)
\end{align*}
from \Cref{Darstelung von Yb}.
\begin{itemize}
\item[(a)] \textsl{Directed edges:} A direct application of \Cref{Erste Charakterisierung gerichtete Kante für MCARMA} gives that $a \rarrow b \notin E_{OG}$, if and only if, neither $Y_a(t)$, $D^{(1)} Y_a(t)$,\ldots, $D^{(p-q-1)}Y_a(t)$ 
nor the integral over the past have any influence on $Y_b(t+h)$. In the representation of the $b$-th component $Y_b(t+h)$, the $a$-th component always vanishes because its coefficient functions are zero. This observation is also evident in \Cref{example} as derived in \Cref{Example Darstellung}.
\item[(b)] \textsl{Undirected edges:} \Cref{Einfache Charakterisierung ungerichtete Kante für MCARMA} yields
\begin{align*}
a \inst b \notin E_{OG}
\:\:\: \Leftrightarrow \:\:\:
\BE[e_a^\top \Oeps(t,h) e_b^\top \Oeps(t, \widetilde{h})]=
\BE[e_a^\top \Oeps(0,h) e_b^\top \Oeps(0, \widetilde{h})]=0, \:\:\: 0\leq h, \widetilde{h}\leq 1.
\end{align*}
This means that the noise terms $e_a^\top  \Oeps(t,h)$ and $e_b^\top  \Oeps(t, \widetilde{h})$ of $Y_a(t+h)$ and $Y_b(t+\widetilde h)$, respectively, are uncorrelated for any $t\geq 0$ and $0\leq h, \widetilde{h}\leq 1$. Again, this observation is also evident in \Cref{example} (\Cref{Example Darstellung}). However, due to the complexity of the expression $\BE[ \varepsilon(t,h) \varepsilon(t,\tilde{h})^\top ] = \int_0^{\min(h,\widetilde{h})} \OBC e^{\BA (h-s)} \BB \BS_L \BB^\top e^{\BA^\top(\widetilde{h}-s)} \OBC^\top ds$, we do not specify the latter. 
\end{itemize}
\end{interpretation}

\begin{interpretation}[Local orthogonality graph]
The interpretation of the directed and the undirected edges in the local orthogonality graph $G_{OG}^0$ is a lot more intricate since the mean square limit of the difference quotient does not exist by definition and \Cref{Remark highes derivative}, respectively, but the limit of the projections does. Therefore we again use the representation for $b\in V$ of \Cref{Lemma 4.7},
\begin{align*}
\texttt{D}_h^{(p-q-1)}Y_b(t,h)
\coloneqq & \frac{D^{(p-q-1)} Y_b(t+h)- D^{(p-q-1)} Y_b(t)}{h} \\
=   & \int_{-\infty}^t  e_b^\top \BM'(0) e^{\BLAM (t-u)} \BT Y_V(u) du 
+ \sum_{m=0}^{p-q-1} e_b^\top \BMnu'(0) \BT D^{(m)} Y_V(t) \\ 
& + e_b^\top O(h)\bR_1+ e_b^\top O(h)\bR_2  + \frac{e_b^\top\Ueps(t,h)}{h}, 
\intertext{and hence,}
P_{\mathcal{L}_{Y_V}(t)} \texttt{D}_h^{(p-q-1)}Y_b(t,h)
=   & \int_{-\infty}^t  e_b^\top \BM'(0) e^{\BLAM (t-u)} \BT Y_V(u) du 
+ \sum_{m=0}^{p-q-1} e_b^\top \BMnu'(0) \BT D^{(m)} Y_V(t) \\ 
& + e_b^\top O(h)\bR_1+ e_b^\top O(h)\bR_2.
\end{align*}
Despite the fact that the $L^2$-limit of $\texttt{D}_h^{(p-q-1)}Y_b(t,h)$ does not exist, the $L^2$-limits of $\sqrt{h}\texttt{D}_h^{(p-q-1)}Y_b(t,h)$  and $P_{\mathcal{L}_{Y_V}(t)} \texttt{D}_h^{(p-q-1)}Y_b(t,h)$ exist.
\begin{itemize}
    \item [(a)]  \textsl{Directed edges:} By \Cref{Erste Charakterisierung gerichtete Kante für MCARMA} we receive that $a \rarrow b \notin E_{OG}^0$, if and only if, neither $Y_a(t)$, $D^{(1)}Y_a(t)$,\ldots, $D^{(p-q-1)}Y_a(t)$ nor the integral over the past have any influence on 
    $\texttt{D}_h^{(p-q-1)}Y_b(t,h)$  if $h$ is small. The same holds for $P_{\mathcal{L}_{Y_V}(t)} \texttt{D}_h^{(p-q-1)}Y_b(t,h)$. \textsl{Given} $\mathcal{L}_{Y_V}(t)$, the $a$-th component does not influence the $b$-th component in the limit, because the corresponding coefficients are zero. Note that in \Cref{example}, we have
    \begin{align*}
        &\int_{-\infty}^t  \BM'(0) e^{\BLAM (t-u)} \BT Y_V(u) du 
        + \BMo'(0) \BT Y_V(t) \\
        & \quad = \int_{-\infty}^t - e^{-(t - u)} 
            \begin{pmatrix*}[c] 
            Y_1(u) \phantom{\: \:- \frac{1}{2} (Y_1(u) - Y_2(u))}\\
            Y_2(u) \phantom{\: \:- \frac{1}{2} (Y_1(u) - Y_2(u))}\\
            Y_3(u) - \frac{1}{2} (Y_1(u) - Y_2(u))
            \end{pmatrix*} du
            + \begin{pmatrix} 
            0\\
            0\\
            0
            \end{pmatrix},
    \end{align*}
    which explains the directed edges in the local orthogonality graph in \Cref{fig: Two graphical models b}.
    
    \item [(b)] \textsl{Undirected edges:} By \Cref{Einfache Charakterisierung ungerichtete Kante für MCARMA} we receive that $a\inst b \notin E_{OG}^0$, if and only if,
    \begin{align*}
         & h \, \BE \left[ \left(\texttt{D}_h^{(p-q-1)}Y_a(t,h) - P_{\mathcal{L}_V(t)}\texttt{D}_h^{(p-q-1)}Y_a(t,h) \right) \right.\\
         & \quad \quad \times 
         \left.\left(\texttt{D}_h^{(p-q-1)}Y_b(t,h) - P_{\mathcal{L}_V(t)}\texttt{D}_h^{(p-q-1)}Y_b(t,h)  \right) \right]\\
         & \quad 
         = \frac{1}{h} \BE \left[ e_a^\top\Ueps(t,h)e_b^\top\Ueps(t,h) \right]
        \stackrel{h\downarrow 0}{\rightarrow} e_a^\top \UBC \BB \Sigma_L \BB^\top \UBC^\top e_b
    \end{align*}
    is zero. Hence, \textsl{given} $\mathcal{L}_{Y}(t)$, $\sqrt{h}\texttt{D}_h^{(p-q-1)}Y_a(t,h)$ and $\sqrt{h} \texttt{D}_h^{(p-q-1)}Y_b(t,h)$ are uncorrelated in the limit. Equivalently, the noise terms $e_a^\top \Ueps(t,h) / \sqrt{h}$ and $e_b^\top \Ueps(t,h) / \sqrt{h}$ are uncorrelated in the limit. Note that in \Cref{example}, we have
    \setlength{\arraycolsep}{4pt}
    \begin{align*}
    \UBC \BB \Sigma_L \BB^\top \UBC^\top
    = \begin{pmatrix}
        2 & 0 & 0 \\
        0 & 2 & 2 \\
        0 & 2 & 3
    \end{pmatrix},
    \end{align*}
    which explains the undirected edges in the local orthogonality graph in \Cref{fig: Two graphical models b}.
\end{itemize}
\end{interpretation}

\begin{remark}\label{remark: comparison to MCAR}
We establish the relationship between our results for ICCSS processes and the results for MCAR processes in \cite{VF23pre}.
\begin{itemize}
    \item[(a)] Since the \textsl{undirected edges} are characterised only by the noise terms $\Oeps(t,h)$ and $\Ueps(t,h)$ and thus, have nothing to do with the inversion of the process, it is not surprising that the characterisations for the undirected edges of the ICCSS$(p,q)$ processes and for the undirected edges of the MCAR$(p)$ processes coincide.
    \item[(b)] In the characterisations of the \textsl{directed edges} of the ICCSS$(p,q)$ process, the case $q=0$ cannot simply be inserted because several matrices become zero-dimensional. However, if we interpret $\OBM(h) e^{\BLAM (t-u)} \BT \: \widehat{=} \: 0_{k\times k}$ if $u<t$, and $\OBMnu(h) \BT \: \widehat{=} \: \BFC e^{\BA h}\BFE_{m+1}$ for $m=0,\ldots,p-1$, as in \Cref{remark: comparison to MCAR 1}, the characterisations of the directed edges for MCAR$(p)$ processes can be seen as special case of \Cref{Zweite Charakterisierung gerichtete Kante für MCARMA} and \Cref{Erste Charakterisierung gerichtete Kante für MCARMA}.
\end{itemize}
\end{remark}

\section{Proofs}\label{sec:proofs}  \label{Appendix:Proofs}
\subsection{Proofs of Subsection~\ref{subsec: state space models}}\label{subsec: Proofs for preliminaries}
\begin{proof}[Proof of \Cref{lemma: uniqueness}]
Assume that there exist two coprime right polynomial fraction descriptions of $H(z)$ as in \eqref{def PQ}, so that
\begin{align*}
    Q_1(z)P_1(z)^{-1}=H(z) = Q_2(z)P_2(z)^{-1}.
\end{align*}
 Then, due to the coprimeness, there exists a matrix polynomial $U(z)$, where $\det(U(z))$ is a non-zero real number \cite[Theorem 16.10]{RU96}, such that 
\begin{align}\label{eq: p=pu}
     P_1(z) = P_2(z) U(z).
\end{align}
Both $P_1(z)$ and $P_2(z)$ have the highest power $I_kz^p$, so $U(z)=I_k$. 
Hence $P_1(z) = P_2(z)$ and finally, $Q_1(z) = Q_2(z)$, which results in the uniqueness of the decomposition. 

The fact that $ H(z)$ is equal to $\BFC \left(z I_{kp} - \BA \right)^{-1} \BB$ follows from the proof of Theorem 3.2 in \cite{BR12}.

Furthermore, the realisations $(\BA, \BB, \BFC)$ and $(\Astar, \Bstar, \Cstar)$ are minimal because $P(z)$ and $Q(z)$ are right coprime and $\deg(\det(P(z))=kp$, see Theorem 6.5-1 of \cite{KA80}. 
Then a consequence of Theorem 2.3.4 in \cite{HA12} is that there exists a non-singular matrix $T$ such that
\begin{align*}
    \BA=T\Astar T^{-1}, \quad \BB=T \Bstar, \quad \text{ and } \quad \BFC=\Cstar T^{-1},
\end{align*}
and for $s<t$,
  \begin{align*} 
        Y(t)&=\Cstar e^{\Astar (t-s)}X(s)+ \int_s^t \Cstar  e^{\Astar (t-u)} \Bstar dL(u)\\
        &=\BFC e^{\BA (t-s)}(TX(s))+ \int_s^t \BFC e^{\BA (t-u)} \BB dL(u).
    \end{align*}
Thus, $\CY$ is a solution of the state space model  $(\Astar, \Bstar, \Cstar, L)$, if and only if, it is a solution of $(\BA, \BB, \BFC, L)$. Finally, $\sigma(\BA)=\sigma(T\Astar T^{-1})=\sigma(\Astar)$.
\end{proof}

\begin{proof}[Proof of \Cref{Lemma: Zusammenhang MCARMA}]
The uniqueness follows directly from \Cref{lemma: uniqueness}. Furthermore, $\pstar=p$ holds by Lemma 6.5-6 of \cite{KA80}.
Since $\Pstar(z)^{-1} \Qstar(z)=Q(z)P(z)^{-1}$ we have $\pstar-\qstar=p-q$ and therefore $\qstar=q$. Comparing the highest order coefficient in $\Qstar(z)P(z)= \Pstar(z) Q(z)$ gives $\Qstar_0=C_q$.  Finally, Lemma 6.3-8 in \cite{KA80} states that $\mathcal{N}(\Pstar)=\mathcal{N}({P})$ and $\mathcal{N}(\Qstar)=\mathcal{N}({Q})$.
\end{proof}

\subsection{Proofs of Section~\ref{sec: Linear prediction of MCARMA processes}}\label{subsec: Proofs for Linear prediction}


\begin{proof}[Proof of \Cref{Xq als integral}]
The proof is divided into four steps. In the first three steps, we derive some auxiliary results which lead in Step 4 to the proof of the statement.\\
\textsl{Step 1:} First, we prove  for all $\varepsilon>0$ and $v\in V$ the asymptotic behaviour
\begin{align}\label{step 1}
\lim_{\vert u \vert \rightarrow \infty} e^{- \varepsilon \vert u \vert} \vert Y_v(u) \vert =0 \quad \text{$\mathbb{P}$-a.s.}
\end{align}
Thus, we relate \eqref{step 1} back to \cite{BR15}, Proposition 2.6, who prove this convergence for stationary univariate CARMA processes that are driven by univariate Lévy processes and whose AR polynomial has no zeros on the \mbox{imaginary axis.}
Therefore, let $\varepsilon>0$ and $v\in V$.  Note that for $t\in \R$, 
\begin{align*}
Y_v(t) 
= \int_{-\infty}^t e_v^\top \OBC e^{\BA (t-u)} \BB dL(u)
= \sum_{\ell=1}^k \int_{-\infty}^t e_v^\top \OBC e^{\BA(t-u)} \BB e_\ell dL_\ell(u)
= \sum_{\ell=1}^k Y_v^\ell(t).
\end{align*}
The process $\CY_v^\ell=(Y_v^\ell (t))_{t\in \R}$ is the stationary solution of the state space model
\begin{align*}
dX(t) = \BA X(t) dt + \BB e_\ell dL_\ell(t), \quad \quad
Y_v^\ell(t) = e_v^\top \OBC X(t),
\end{align*}
and has the transfer function 
\begin{align*}
H_v^\ell(z)=e_v^\top \OBC (zI_{kp} -\BA)^{-1} \BB e_{\ell}. 
\end{align*}
Then \cite{KA80} provides in Lemma 6.3-8 the existence of (right) coprime  polynomials $P_v^\ell(z)$ and $Q_v^\ell(z)$ (polynomials with no common zeros) as in \eqref{def PQ} so that $H_v^\ell(z) = Q_v^\ell(z)/P_v^\ell(z)$. Note that in the univariate setting the problem of the existence of a coprime right polynomial fraction description of the form \eqref{def PQ} does not arise. Indeed, here $1 \cdot p=\deg(\det(P_v^\ell(z))=\deg(P_v^\ell(z))$ follows immediately, and the constant before the $p$-th power can be included in $Q_v^\ell(z)$ without loss of generality so that $P_v^\ell(z)$ is a polynomial of degree $p$ that has a $1$ as the leading coefficient. Thus, the classes of univariate CARMA processes and univariate causal continuous-time state space models are equivalent \cite[Corollary 3.4]{SC12} implying that  $\CY_v^\ell$ is a univariate  CARMA process driven by a univariate Lévy process. 
 Now, \cite{BE09}, Definition 4.7.1, provides that the poles of $H_v^\ell(z)$ are the roots of $P_v^\ell(z)$ including multiplicity. In addition, \cite{BE09}, Theorem 12.9.16, yields that the poles of $H_v^\ell(z)$ are a subset of $\sigma(\BA)$ resulting in $$\mathcal{N}(P_v^\ell) = \{z\in \C: P_v^\ell(z)=0\}  \subseteq \sigma(\BA ) \subseteq (-\infty,0) +i\R,$$ which means that the AR polynomial $P_v^\ell(z)$  has  no zeros on the imaginary axis. 
Thus, $\CY_v^\ell$ satisfies the assumptions in  \cite{BR15}, Proposition 2.6, and we obtain for $\ell=1,\ldots,k$ that
\begin{align*}
\lim_{\vert u \vert \rightarrow \infty} e^{- \varepsilon \vert u \vert} Y_v^\ell(u) =0 \quad \text{$\mathbb{P}$-a.s.}
\end{align*}
Therefore,
\begin{align*}
\lim_{\vert u \vert \rightarrow \infty} e^{- \varepsilon \vert u \vert} Y_v(u) =
\sum_{\ell=1}^k \lim_{\vert u \vert \rightarrow \infty} e^{- \varepsilon \vert u \vert} Y_v^\ell(u) =0 \quad \text{$\mathbb{P}$-a.s.},
\end{align*}
and finally, the claim \eqref{step 1} follows.\\
\textsl{Step 2:} Next, we show that 
\begin{align}\label{step 2}
\lim_{s\rightarrow -\infty} \int_{s}^t  e^{-\lambda (t-u)} \vert Y_v(u) \vert du
\end{align}
exists $\mathbb{P}$-a.s.~for $t\in \R$ and $\lambda>0$. 

From \eqref{step 1} we obtain that there exists some set $\Omega_0 \in \mathcal{F}$ with $\mathbb{P}(\Omega_0) =1$ such that for all $\omega \in \Omega_0$ and $\gamma >0$ there exists a $u_0(\omega)<0$ with
\begin{align*}
e^{\frac{\lambda}{2} u} \vert  Y_v(\omega, u) \vert = e^{- \frac{\lambda}{2} \vert u \vert} \vert  Y_v(\omega, u) \vert  \leq \gamma \quad \forall \: u \leq u_0(\omega).
\end{align*}
Then we obtain for $s< u_0(\omega)$ that
\begin{align*}
\int_s^t e^{-\lambda(t-u)} \vert Y_v(\omega, u) \vert du 
&= \int_{u_0(\omega)}^t e^{-\lambda(t-u)} \vert Y_v(\omega, u) \vert du  
+ \int_s^{u_0(\omega)} e^{-\lambda(t-u)} \vert Y_v(\omega, u) \vert du  \\
& \leq \int_{u_0(\omega)}^t e^{-\lambda(t-u)} \vert Y_v(\omega,u) \vert du  + \gamma e^{-\lambda t} \frac{2}{\lambda}.
\end{align*}
Thus, by dominated convergence the limit in \eqref{step 2}
exists $\mathbb{P}$-a.s.~for $t\in \R$ and $\lambda>0$.\\

\textsl{Step 3:} Eventually, we derive that not only the univariate integral \eqref{step 2} exists, but also 
\begin{align}\label{step 3}
\lim_{s\rightarrow -\infty} \int_{s}^t e^{\BLAM (t-u)} \BT Y(u) du 
\end{align}
exists $\mathbb{P}$-a.s.~for $t\in \R$. First, Assumption \eqref{Assumption Q} provides that $\sigma(\BLAM ) \subseteq (-\infty,0) +i\R$ and thus, 
$\text{spabs}(\BLAM) \coloneqq \max\{ \Re(\lambda) : \lambda \in \sigma(\BLAM ) \} <0,
$
where $\Re(\lambda)$ denotes the real part of $\lambda$. Therefore, there exists a $-\lambda \in (\text{spabs}(\BLAM ),0)$. Then \cite{BE09}, Proposition 11.18.8, provides a constant $c_1>0$ such that
\begin{align}\label{Abschätzung von exp}
\Vert e^{\BLAM t} \Vert \leq c_1 e^{-\lambda t} \quad \forall \: t\geq 0.
\end{align}
Now, we obtain
\begin{align*}
\left\| \int_s^t e^{\BLAM (t-u)} \BT Y(u) du \right\|
\leq  c_1 \Vert \BT \Vert  \sum_{v\in V}  \int_s^t e^{-\lambda (t-u)} \vert Y_v(u)\vert du.
\end{align*}
Due to \eqref{step 2} the limit of each of those addends exists, so \eqref{step 3} exists $\mathbb{P}$-a.s.~for $t\in \R$. \\
\textsl{Step 4:} Finally, we are able to prove the statement of the proposition. Recall that due to \eqref{erste komponenten} for $s,t \in \R$, $s<t$,
\begin{align*}
X^q(t)=e^{\BLAM (t-s)}X^q(s)+\int_s^t e^{\BLAM (t-u)} \BT Y(u)du.
\end{align*}
 Since we assume that $\CX$ is the unique stationary solution of the stochastic differential equation  \eqref{controller dgl},  $\CX^q$ is also strictly stationary and $X^q(s)$ and $X^q(0)$ have the same distribution for all $s \in \R$. Moreover, it follows from Assumption \eqref{Assumption Q} that $\sigma(\BLAM ) \subseteq (-\infty,0) +i\R$. These properties lead to 
\begin{align*}
e^{\BLAM (t-s)} X^q(s) \rightarrow 0_{kq} \quad  \text{ as } s \rightarrow - \infty,
\end{align*}
in distribution and in probability by Slutsky's lemma, since the limit is a degenerate random vector. In combination with \eqref{step 3} we receive for  $t\in \R$ the statement
\begin{equation*}
\lim_{s\rightarrow -\infty} \left(e^{\BLAM (t-s)} X^q(s)+ \int_s^t e^{\BLAM (t-u)} \BT Y (u)du\right)
= \int_{-\infty}^t e^{\BLAM (t-u)} \BT Y(u)du \quad \mathbb{P}\text{-a.s.} \qedhere
\end{equation*}
\end{proof}

\begin{proof}[Proof of \Cref{Darstelung von Yb}]
Let $t\in \R$, $h\geq 0$, and $a \in V$. First of all, due to \eqref{eq: def CC}, we receive
\begin{align*}
Y_a(t+h) 
= e_a^\top  \OBC X(t+h) \quad \text{ and } \quad
D^{(p-q-1)}Y_a(t+h)
= e_a^\top  \UBC X(t+h).
\end{align*}
From now on, the proofs for the two representations differ only in the choice of $\OBC$ and $\UBC$, respectively. Therefore, we will only continue with the representation of $Y_a(t+h)$. Due to 
\eqref{X aufgeteilt} we have
\begin{align*}
  Y_a(t+h) 
  =   e_a^\top  \OBC \left( e^{\BA h} X(t) + \int_t^{t+h} e^{\BA (t+h-u)} \BB dL(u) \right) 
  =  e_a^\top  \OBC e^{\BA h} X(t) +  e_a^\top  \Oeps(t,h).
\end{align*}
Here,
\begin{align*}
e_a^\top  \OBC e^{\BA h} X(t) 
&=  e_a^\top  \OBC e^{\BA h} \begin{pmatrix}
    X^q (t), & X^{(q+1)}(t), & \cdots, & X^{(p)}(t)
\end{pmatrix}^\top \\
&= e_a^\top  \OBC e^{\BA h} \OBE  X^q (t) + \sum_{j=1}^{p-q} e_a^\top  \OBC e^{\BA h} \BFE_{q+j}  X^{(q+j)} (t).
\end{align*}
\Cref{restliche komponenten} and interchanging the summation order
imply

\begin{align} \label{eqa1}
&e_a^\top  \OBC e^{\BA h} X(t) \nonumber \\
&\quad = e_a^\top  \OBC e^{\BA h} \OBE  X^q (t) + \sum_{j=1}^{p-q} e_a^\top  \OBC e^{\BA h} \BFE_{q+j} \UBE^\top  \left( \BLAM^j X^q(t) + \sum_{m=0}^{j-1}\BLAM^{j-1-m} \BT D^{(m)} Y(t)\right)  \nonumber \\
&\quad= e_a^\top  \OBM(h) X^q (t) + \sum_{m=0}^{p-q-1} e_a^\top  \OBMnu (h) \BT D^{(m)} Y(t).
\end{align}
Finally, we obtain due to 
\Cref{Xq als integral}, 
\begin{equation*}
e_a^\top  \OBC e^{\BA h} X(t) 
= \int_{-\infty}^t e_a^\top  \OBM(h) e^{\BLAM (t-u)} \BT Y(u) du  + \sum_{m=0}^{p-q-1} e_a^\top  \OBMnu (h) \BT D^{(m)} Y(t) \quad \text{$\mathbb{P}$-a.s.}  \qedhere
\end{equation*}
\end{proof}

\begin{proof}[Proof of \Cref{Xq als L2 integral}]
Let $a,v\in V$ and define $F(t) = e_a^\top e^{\BLAM t} \BT e_v$
for $t\geq 0$. First, for ${s},t\in \R$, ${s}<t$,
\begin{align*}
\lim_{n\rightarrow \infty} \frac{t-s}{n} \sum_{\ell=1}^n F\left(t-{{s}}-\ell \: \frac{t-{s}}{n}\right) Y_v \left( {s}+\ell \: \frac{t-{s}}{n} \right) = \int_{{s}}^t F(t-u) Y_v(u) du \quad \text{$\mathbb{P}$-a.s.}
\end{align*}
 due to the definition of the integral. Using the theorem of dominated convergence, we show that this convergence also holds in the $L^2$ sense. Indeed, from  the triangle inequality 
\begin{align*}
&\: \left| \int_{s}^t F(t-u) Y_v(u) du -  \frac{t-s}{n} \sum_{\ell=1}^n F\left(t-{s}-\ell \: \frac{t-{s}}{n}\right) Y_v \left( {s}+\ell \: \frac{t-{s}}{n} \right) \right| \\
& \quad\quad\leq \:  \int_{s}^t \vert F(t-u) \vert \vert Y_v(u) \vert du +  \frac{t-s}{n} \sum_{\ell=1}^n \left| F\left(t-{s}-\ell \: \frac{t-{s}}{n}\right)\right| \: \left| Y_v\left({s}+\ell \: \frac{t-{s}}{n} \right) \right| \\
& \quad\quad \leq \: 2 (t-{s}) \left( \sup_{u \in [0, t-{s}]}  \vert F(u) \vert \right) \left(\sup_{u\in [{s},t]} \vert Y_v(u) \vert \right)
\end{align*}
follows.
This majorant is integrable, because 
\begin{align*}
\sup_{u\in [0,t-{s}]} \vert Y_v(u) \vert 
= \sup_{u\in [0,t-{s}]} \vert e_v^\top  \OBC X(u) \vert
\leq \sup_{u\in [0,t-{s}]} \Vert e_v^\top  \OBC \Vert \Vert X(u) \Vert
\leq c \sup_{u\in [0,t-{s}]} \Vert X(u) \Vert,
\end{align*}
for some constant $c\geq0$ and thus,

\begin{align} \label{eq21}
\BE \left[ \left(\sup_{u\in [0,t-{s}]} \vert Y_v(u) \vert\right)^2 \right]
\leq c^2 \:  \BE \left[ \left( \sup_{u\in [0,t-{s}]} \Vert X(u) \Vert \right)^2 \right] < \infty,
\end{align}
due to \Cref{Assumption on Levy process} and \cite{BR12}, Lemma A.4. Furthermore,  $\sup_{u\in [0, t-{s}]} \vert F(u) \vert<\infty$ since $F$ is a continuous function. In summary,
\begin{align*}
\int_{s}^t F(t-u) Y_v(u) du 
= \limm \frac{t-s}{n} \sum_{\ell=1}^n F\left(t-{s}-\ell \frac{t-{s}}{n}\right) Y_v \left( {s}+\ell \frac{t-{s}}{n} \right).
\end{align*}
For the second step of this proof, we recall that for $t\in \R$,
\begin{align*}
\int_{-\infty}^t F(t-u) Y_v(u) du 
= \lim_{{s}\rightarrow -\infty} \int_{{s}}^t F(t-u) Y_v(u) du \quad \text{$\mathbb{P}$-a.s.},
\end{align*}
due to \Cref{Xq als integral}.
 Again, using the theorem of dominated convergence, we show that this convergence holds in the $L^2$ sense. For ${s}<t$ it follows that 
 
\begin{align*}
  \: \left| \int_{-\infty}^t F(t-u) Y_v(u) du - \int_{{s}}^t F(t-u) Y_v(u) du  \right| 
&  \leq \: \int_{t-{s}}^{\infty} \vert  F(u) \vert  \vert  Y_v(t-u) \vert du \\ 
& \leq \sum_{n=0}^{\infty} \sup_{u \in [n,n+1]} \left| F(u) \right|  \sup_{u \in [n,n+1]} \vert  Y_v(t-u) \vert.
\end{align*}
To see that this majorant is in $L^2$, we use Fubini, Cauchy-Schwarz inequality and the stationarity of $\CY$. This yields to
\begin{align*}
&\: \BE \left[ \left( \sum_{n=0}^{\infty} \sup_{u \in [n,n+1]} \vert F(u) \vert  \sup_{u\in [n,n+1]} \vert  Y_v(t-u) \vert \right)^2 \right] \\
& \quad\quad\leq  \sum_{n=0}^{\infty} \sum_{m=0}^{\infty} \sup_{u \in [n,n+1]} \vert F(u) \vert \sup_{u \in [m,m+1]} \vert F(u) \vert  \\
& \quad\quad \quad\quad \times \left( \BE \left[ \left( \sup_{u\in [n,n+1]} \vert  Y_v(t-u) \vert \right)^2 \right]
\BE \left[ \left( \sup_{u\in [m,m+1]} \vert  Y_v(t-u) \vert \right)^2 \right] \right)^{1/2} \\
& \quad\quad = \: \left(  \sum_{n=0}^{\infty} \sup_{u \in [n,n+1]} \vert F(u) \vert \right)^2 \BE \left[ \left(  \sup_{u\in [0,1]} \vert  Y_v(u) \vert \right)^2 \right] < \infty,
\end{align*}
where we used \eqref{eq21} and $\sum_{n=0}^{\infty} \sup_{u \in [n,n+1]} \vert F(u) \vert < \infty$ by the definition of $F$ and \eqref{Abschätzung von exp}. In summary, we obtain
\begin{equation*}
e_a^\top  \int_{-\infty}^t  e^{\BLAM (t-u)} \BT e_v Y_v(u) du
= \int_{-\infty}^t F(t-u) Y_v(u) du 
= \underset{{s}\rightarrow -\infty}{\text{l.i.m.\:}} \int_{{s}}^t F(t-u) Y_v(u) du,
\end{equation*}
and the integral is in $\mathcal{L}_{Y_v}(t)$. The existence of $X^q(t)$ as an $L^2$-limit follows immediately from this.
\end{proof}

\begin{proof}[Proof of \Cref{Lemma 4.7}]
Recall that due to \Cref{Darstelung von Yb} and $\Ueps(t,0)=0_k \in \R^k$
\begin{align}
&\frac{D^{(p-q-1)} Y_a(t+h)- D^{(p-q-1)} Y_a(t)}{h} \nonumber \\
&\quad = \int_{-\infty}^t  e_a^\top \frac{\BM(h)-\BM(0)}{h} e^{\BLAM (t-u)} \BT Y(u) du \nonumber \\
&\quad \phantom{=} + \sum_{m=0}^{p-q-1} e_a^\top \frac{\BMnu(h)-\BMnu(0)}{h} \BT D^{(m)} Y(t) +  e_a^\top  \frac{\Ueps(t,h)}{h} \quad \text{$\mathbb{P}$-a.s.} \label{eqac}
\end{align}
Replacing the matrix exponential with its power series, it holds that
\begin{align} \label{eqab}
\begin{split}
    \frac{\BM(h)-\BM(0)}{h}
    &= \UBC \frac{e^{\BA h}-I_{kp}}{h} \Bigg( \OBE + \sum_{j=1}^{p-q} \BFE_{q+j} \UBE^\top  \BLAM^j \Bigg) \\
    &=\BM'(0)+O(h) \Bigg( \OBE + \sum_{j=1}^{p-q} \BFE_{q+j} \UBE^\top  \BLAM^j \Bigg), \\
    \frac{\BMnu(h)-\BMnu(0)}{h}    
    &= \UBC \frac{e^{\BA h}-I_{kp}}{h} \sum_{j=m + 1}^{p-q} \BFE_{q+j} \UBE^\top  \BLAM^{j-1-m} \\
    &=\BMnu'(0) + O(h)  \sum_{j=m + 1}^{p-q} \BFE_{q+j} \UBE^\top  \BLAM^{j-1-m}.
\end{split}
\end{align}
Furthermore, we define
\begin{align}
\begin{split}\label{eqhilf12}
\bR_1 &= \int_{-\infty}^t \Bigg( \OBE + \sum_{j=1}^{p-q} \BFE_{q+j} \UBE^\top  \BLAM^j \Bigg) e^{\BLAM (t-u)} \BT Y(u) du,\\
\bR_2 &= \sum_{m=0}^{p-q-1} \sum_{j=m + 1}^{p-q} \BFE_{q+j} \UBE^\top  \BLAM^{j-1-m} \BT D^{(m)} Y(t).
\end{split}
\end{align}
If we plug \eqref{eqab} and \eqref{eqhilf12} in \eqref{eqac} we obtain the stated representation. Moreover, from \Cref{Xq als L2 integral}  we know that $\bR_1$ is in $\mathcal{L}_{Y}(t)$ and from \Cref{Remark Derivatives in linear space} we receive that $\bR_2$ is in $\mathcal{L}_{Y}(t)$. Since  $\mathcal{L}_{Y(t)}$  and $(L(s)-L(t))_{t \leq s \leq t+h}$ are independent we receive that $R_1,R_2\in\mathcal{L}_{Y}(t)$ are independent of $\Ueps(t,h)$. Finally, 
\begin{eqnarray*}
    \frac{1}{h}\BE\left[ (e_a^\top \Ueps(t,h))^2 \right] = \frac{1}{h} e_a^\top \UBC   \int_0^{h} e^{\BA u} \BB \Sigma_L \BB^\top  e^{\BA^\top  u} du \UBC e_a 
\stackrel{h\downarrow 0}{\longrightarrow} e_a^\top  \UBC \BB \Sigma_L \BB^\top  \UBC e_a.
\end{eqnarray*}
$\UBC \BB \Sigma_L \BB^\top  \UBC$ is positive definite due to $\Sigma_L>0$ and $\UBC$, $\BB$ being of full rank by Assumption \eqref{Assumption Q}. Therefore, the limit $e_a^\top \UBC \BB \Sigma_L \BB^\top  \UBC e_a>0$  and, of course, $\BE[(e_a^\top \Ueps(t,h))^2]/h^2 $ converges then to infinity. 
\end{proof}

\begin{proof}[Proof of \Cref{Projektionen für MCARMA}]
Based on \Cref{Darstelung von Yb}, the proofs of the two orthogonal projections differ only in the choice of 
$\OBM(\cdot)$ or $\BM(\cdot)$,  $\OBMnu(\cdot)$ or $\BMnu(\cdot)$, and  $ \Oeps(\cdot,\cdot)$ or $ \Ueps(\cdot,\cdot)$. Thus, we only prove the representation of $P_{\mathcal{L}_{Y_S}(t)} Y_a(t+h)$. 
Let $h\geq 0$, $t\in \R$, $S \subseteq V$, and $a\in V$. From \Cref{Darstelung von Yb} recall that $\mathbb{P}$-a.s.
\begin{align*}
Y_a(t+h)
= & \: \int_{-\infty}^t e_a^\top  \OBM(h) e^{\BLAM (t-u)} \BT Y(u) du + \sum_{m=0}^{p-q-1} e_a^\top  \OBMnu (h) \BT D^{(m)} Y(t) +e_a^\top  \Oeps(t,h).
\end{align*}
We calculate the projections of the summands separately. For the first \mbox{summand we get}
\begin{align*}
&P_{\mathcal{L}_{Y_S}(t)} \left( \int_{-\infty}^t e_a^\top  \OBM(h) e^{\BLAM (t-u)} \BT Y(u) du \right) \\
&\quad=\sum_{v\in S} \int_{-\infty}^t e_a^\top  \OBM(h) e^{\BLAM(t-u)} \BT e_v Y_v(u) du \\
&\quad \quad\phantom{h} + P_{\mathcal{L}_{Y_S}(t)} \left( \sum_{v\in V\setminus S} \int_{-\infty}^t e_a^\top  \OBM(h) e^{\BLAM (t-u)} \BT e_v Y_v(u) du \right),
\end{align*}
since, because of \Cref{Xq als L2 integral}, the integrals are in $\mathcal{L}_{Y_S}(t)$ for $v\in S$. For the second summand, we obtain
\begin{align*}
&P_{\mathcal{L}_{Y_S}(t)} \left( \sum_{m=0}^{p-q-1} e_a^\top  \OBMnu (h) \BT D^{(m)} Y(t) \right) \\
& \quad = \sum_{v\in S} \sum_{m=0}^{p-q-1}  e_a^\top  \OBMnu (h) \BT e_v D^{(m)} Y_v(t)\\
& \phantom{\quad \quad}  
+ P_{\mathcal{L}_{Y_S}(t)} \left(  \sum_{v\in V\setminus S} \sum_{m=0}^{p-q-1} e_a^\top  \OBMnu (h) \BT e_v D^{(m)} Y_v(t) \right),
\end{align*}
since, due to \Cref{Remark Derivatives in linear space}, the derivatives of $Y_v(t)$ for $v\in S$ are in $\mathcal{L}_{Y_S}(t)$.
For the third summand $e_a^\top  \Oeps(t,h)$ we note that $(Y_S(s))_{s\leq t}$ and $(L(s)-L(t))_{t \leq s \leq t+h}$ are independent. We obtain immediately that
$P_{\mathcal{L}_{Y_S}(t)} e_a^\top  \Oeps(t,h) =0$. If we put all three summands together, we get the assertion. 
\end{proof}

\begin{proof}[Proof of \Cref{Projektionen mit limneu}]
Let $S \subseteq V$, $a\in V$, $h\geq 0$, and $t\in \R$. First of all, due to \Cref{Projektionen für MCARMA} and similar ideas as in \eqref{eqa1},
\begin{align*}
&P_{\mathcal{L}_Y(t)} D^{(p-q-1)}Y_a(t+h) \\
&\quad = \int_{-\infty}^t e_a^\top \BM(h)  e^{\BLAM (t-u)} \BT Y(u) du +\sum_{m=0}^{p-q-1} e_a^\top \BMnu (h) \BT D^{(m)} Y(t) \\
&\quad = e_a^\top  \UBC e^{\BA h} X(t).
\end{align*}
Then, due to 
     \begin{align*}
         &\limhh \BE \left[ \left( 
         P_{\mathcal{L}_{Y}(t)} \left(\frac{D^{(p-q-1)} Y_a(t+h)- D^{(p-q-1)} Y_a(t)}{h}\right) -e_a^\top  \UBC \BA X(t)
         \right)^2\right] \\
         &\quad =\limhh \BE \left[ \left( 
         e_a^\top  \UBC \frac{e^{\BA h}-I_{kp}}{h} X(t) -e_a^\top  \UBC \BA X(t)
         \right)^2\right] \\
         &\quad = \limhh e_a^\top  \UBC \left(\frac{e^{\BA h}-I_{kp}}{h}-\BA \right) c_{XX}(0) \left(\frac{e^{\BA h}-I_{kp}}{h}-\BA \right)^\top  \UBC^\top  e_a 
         =0,
     \end{align*}   
we obtain
     \begin{align*}
         \limh P_{\mathcal{L}_{Y}(t)} \left(\frac{D^{(p-q-1)} Y_a(t+h)- D^{(p-q-1)} Y_a(t)}{h}\right)
         = e_a^\top  \UBC \BA X(t)  \quad\text{$\mathbb{P}$-a.s.}
     \end{align*}
Together with \cite{BR91}, Proposition 2.3.2.(iv,vii), it follows that
     \begin{align*}
         &\limh P_{\mathcal{L}_{Y_{S}}(t)} \left(\frac{D^{(p-q-1)} Y_a(t+h)- D^{(p-q-1)} Y_a(t)}{h}\right) \\
         & \quad = \limh P_{\mathcal{L}_{Y_{S}}(t)} P_{\mathcal{L}_{Y}(t)} \left(\frac{D^{(p-q-1)} Y_a(t+h)- D^{(p-q-1)} Y_a(t)}{h}\right)\\
         &\quad = P_{\mathcal{L}_{Y_{S}}(t)} \left( e_a^\top  \UBC \BA X(t)\right)\quad \text{$\mathbb{P}$-a.s.}
     \end{align*}
Again, similar to the proof of \eqref{eqa1},
\begin{equation*}
e_a^\top  \UBC \BA X(t) 
= \int_{-\infty}^t e_a^\top  \BM'(0) e^{\BLAM (t-u)} \BT Y(u) du  + \sum_{m=0}^{p-q-1} e_a^\top  \BMnu' (0) \BT D^{(m)} Y(t) 
\quad \text{$\mathbb{P}$-a.s.}
\end{equation*}
We obtain  replacing $\BM(h)$ by $\BM'(0)$ and  $\BMnu(h)$ by $\BMnu'(0)$ in the proof of \Cref{Projektionen für MCARMA}, 
\begin{align*}
&\limh P_{\mathcal{L}_{Y_{S}}(t)} \left(\frac{D^{(p-q-1)} Y_a(t+h)- D^{(p-q-1)} Y_a(t)}{h}\right) \\
&\quad =\: 
\sum_{v\in S} \int_{-\infty}^t e_a^\top    \BM'(0) e^{\BLAM(t-u)} \BT e_v Y_v(u) du 
+ \sum_{v \in S} \sum_{m=0}^{p-q-1} e_a^\top  \BMnu'(0) \BT e_v D^{(m)} Y_v(t) \\
&\quad\quad\phantom{=} + P_{\mathcal{L}_{Y_{S}}(t)} \left(
\sum_{v\in V\setminus S} \int_{-\infty}^t e_a^\top    
\BM'(0) e^{\BLAM(t-u)} \BT e_v Y_v(u) du 
\right) \\
&\quad\quad\phantom{=} + P_{\mathcal{L}_{Y_{S}}(t)} \left(\sum_{v \in V\setminus S} \sum_{m=0}^{p-q-1} e_a^\top 
\BMnu'(0) \BT e_v D^{(m)} Y_v(t) \right) \quad \mathbb{P}\text{-a.s.} 
\end{align*}
as claimed.\\
The second assertion follows directly from \Cref{Darstelung von Yb} and \Cref{Projektionen für MCARMA}, which give
\begin{align*}
&  D^{(p-q-1)} Y_a(t+h) - P_{\mathcal{L}_{Y}(t)} D^{(p-q-1)}Y_a(t+h) \\
&\quad =  \: \int_{-\infty}^t  e_a^\top  \BM(h) e^{\BLAM (t-u)} \BT Y(u) du + \sum_{m=0}^{p-q-1} e_a^\top  \BMnu(h)^\top  \BT D^{(m)} Y(t) 
+  e_a^\top  \Ueps(t,h)\\
&\quad\quad \phantom{=} - \int_{-\infty}^t e_a^\top \BM(h) e^{\BLAM (t-u)} \BT Y(u) du - \sum_{m=0}^{p-q-1} e_a^\top \BMnu (h) \BT D^{(m)} Y(t) \\
&\quad= e_a^\top  \Ueps(t,h). \qedhere
\end{align*}
\end{proof}

\subsection{Proofs of Section~\ref{sec: orthogonality graph for MCARMA processes}}\label{subsec: Proofs for orthogonality graph}
\begin{proof}[Proof of \Cref{Graph well defined}] $\mbox{}$
\begin{itemize}
    \item[(A.2)] The proof of Assumption \ref{Assumption an Dichte} is elaborate and has already been given in \cite{VF23pre}, Proposition 6.5, for MCAR$(p)$ processes. It can be directly generalised to ICCSS$(p,q)$ processes, so we do not give the full proof. We simply note that we only require that $Q(i\lambda) P(i\lambda)^{-1}$ has full rank and $\BS_L>0$ to obtain that $f_{YY}(\lambda)>0$  for $\lambda \in \R$. 
    Indeed,  Assumption \eqref{Assumption Q} provides that $Q(i\lambda)$ is of full rank and $\mathcal{N}(P)\subseteq (-\infty,0)+i\R$, so we directly receive that $Q(i\lambda) P(i\lambda)^{-1}$ has full rank as well. 
    Furthermore, we require that $\sigma(\BA)\subseteq (-\infty,0)+i\R$, but this is also true due to  Assumption \eqref{Assumption Q}.
    Finally, it is a necessity that  $\OBC \BB \BS_L \BB^\top \OBC^\top >0$.
    Again, $\BS_L>0$, $\OBC$ is of full rank by Assumption \eqref{Assumption Q}, and $\BB$ is of full rank by definition, so $\OBC \BB \BS_L \BB^\top \OBC^\top >0$.

    \item[(A.3)] For Assumption \ref{Assumption purely nondeterministic of full rank} we apply that $\sigma(\BA ) \subseteq (-\infty, 0) + i \R$ and hence
    \begin{align*}
        \limh P_{\mathcal{L}_{X}(t)}X(t+h) = \limh  e^{\BA h} X(t) =0,
    \end{align*}
    resulting in $\CX$ being purely nondeterministic. By \cite{RO67}, III, (2.1) and Theorem 2.1 this is equivalent to $\bigcap_{t\in \R} \mathcal{L}_X(t)=\{0\}$. Since $\bigcap_{t\in \R} \mathcal{L}_{Y}(t) \subseteq \bigcap_{t\in \R} \mathcal{L}_X(t)$ the process $\CY$ is purely nondeterministic as well.
\end{itemize}
 Finally, the Markov properties follow from \cite{VF23pre}, Section~5; see also \cite{VF23pre}, Proposition 6.6 and Proposition 6.7 for MCAR$(p)$ processes.
\end{proof}

Next, we prove \Cref{Erste Charakterisierung gerichtete Kante für MCARMA}, since the proof of \Cref{Zweite Charakterisierung gerichtete Kante für MCARMA} is based on \Cref{Erste Charakterisierung gerichtete Kante für MCARMA}.

\begin{proof}[Proof of \Cref{Erste Charakterisierung gerichtete Kante für MCARMA}] \phantom{a}\\
(a) \, Recall that due to \Cref{Charakterisierung als Gleichheit der linearen Vorhersage} we have no directed edge $a\rarrow b \notin E_{OG}$, if and only if, 
for $0 \leq h \leq 1$ and $t\in \R$,
\begin{align*}
&P_{\mathcal{L}_{Y}(t)}Y_b(t+h) = P_{\mathcal{L}_{Y_{V\setminus\{a\}}}(t)}Y_b(t+h)  \quad \text{$\mathbb{P}$-a.s.}
\end{align*}

From \Cref{Projektionen für MCARMA} we obtain
for $0\leq h \leq 1$ and $t\in \R$,
\begin{align*}
P_{\mathcal{L}_{Y}(t)}Y_b(t+h)
= &\sum_{v \in V} \int_{-\infty}^t e_b^\top  \OBM(h)  e^{\BLAM(t-u)} \BT e_v Y_v(u) du \\
&+ \sum_{v \in V} \sum_{m=0}^{p-q-1} e_b^\top  \OBMnu (h) \BT e_v D^{(m)} Y_v(t),\\
P_{\mathcal{L}_{Y_{{V\setminus\{a\}}}}(t)} Y_b(t+h)
= &\sum_{v \in V\setminus\{a\}} \int_{-\infty}^t e_b^\top  \OBM(h)  e^{\BLAM(t-u)} \BT e_v Y_v(u) du \\
  &\quad + \sum_{v \in V\setminus\{a\}} \sum_{m=0}^{p-q-1} e_b^\top  \OBMnu (h) \BT e_v D^{(m)} Y_v(t) \\
  &\quad + P_{\mathcal{L}_{Y_{V\setminus\{a\}}}(t)} \left( \int_{-\infty}^t e_b^\top  \OBM(h) e^{\BLAM(t-u)} \BT e_a Y_a(u) du \right) \\
  &\quad + P_{\mathcal{L}_{Y_{V\setminus\{a\}}}(t)} \left(\sum_{m=0}^{p-q-1} e_b^\top  \OBMnu (h) \BT e_a D^{(m)} Y_a(t) \right) \quad \text{$\mathbb{P}$-a.s.}
\end{align*}
 We equate 
the two orthogonal projections and remove the coinciding terms. Then we receive that $a \rarrow b \notin E_{OG}$, if and only if, for $0\leq h\leq 1$ and $t\in \R$,
\begin{align*}
& \int_{-\infty}^t e_b^\top  \OBM(h) e^{\BLAM(t-u)} \BT e_a Y_a(u) du + \sum_{m=0}^{p-q-1} e_b^\top  \OBMnu (h) \BT e_a D^{(m)} Y_a(t) \\ 
 & \quad= \: P_{\mathcal{L}_{Y_{V\setminus\{a\}}}(t)} \left( \int_{-\infty}^t e_b^\top  \OBM(h) e^{\BLAM(t-u)} \BT e_a Y_a(u) du \right) \\
& \quad\quad \phantom{h} + P_{\mathcal{L}_{Y_{V\setminus\{a\}}}(t)} \left(\sum_{m=0}^{p-q-1} e_b^\top  \OBMnu (h) \BT e_a D^{(m)} Y_a(t) \right) \quad \text{$\mathbb{P}$-a.s.}
\end{align*}
 The expression on the left hand side is in $\mathcal{L}_{Y_{a}}(t)$ and the expression on the right side is in $\mathcal{L}_{Y_{V\setminus\{a\}}}(t)$. Since $\mathcal{L}_{Y_{V\setminus\{a\}}}(t) \cap \mathcal{L}_{Y_{a}}(t)= \{0\}$ due to \eqref{eq: linear independence}, $a \rarrow b \notin E_{OG}$, if and only if, for $0\leq h\leq 1$ and $t\in \R$,
\begin{align}\label{eq: gegeben für gerichtete Kante}
\int_{-\infty}^t e_b^\top  \OBM(h) e^{\BLAM (t-u)} \BT e_a Y_a(u) du + \sum_{m=0}^{p-q-1} e_b^\top  \OBMnu (h) \BT e_a D^{(m)} Y_a(t) =0 \quad \text{$\mathbb{P}$-a.s.}
\end{align}
In the following, we show that this characterisation is in turn equivalent to
\begin{align}\label{eq: behauptung für gerichtete Kante}
e_b^\top  \OBM(h) e^{\BLAM t} \BT e_a = 0 \quad \text{and} \quad 
e_b^\top \OBMnu (h) \BT e_a = 0,
\end{align}
for $m=0,\ldots,p-q-1$, $0 \leq h \leq 1$, and $t \geq 0$.

If \eqref{eq: behauptung für gerichtete Kante} holds, we immediately obtain that \eqref{eq: gegeben für gerichtete Kante} is valid. 
Now, suppose \eqref{eq: gegeben für gerichtete Kante} holds. We convert the two summands in \eqref{eq: gegeben für gerichtete Kante} into their spectral representation. 
Hence, note that due \cite{BE09}, Proposition 11.2.2, and $\sigma(\BLAM) \subseteq (-\infty, 0) + i\R$ the equality
\begin{align*}
\uint e^{-i \lambda s} {\boldsymbol 1}_{\{s\geq 0\}} e_b^\top  \OBM(h) e^{\BLAM s} \BT e_a ds
&= e_b^\top  \OBM(h) (i\lambda I_{kq} - \BLAM )^{-1}\BT e_a, \quad\lambda\in\R,
\end{align*}
holds.
Now \cite{RO67} I, Example 8.3, provides the spectral representation of the first summand
\begin{align*}
\int_{-\infty}^t e_b^\top  \OBM(h)  e^{\BLAM (t-u)} \BT e_a Y_a(u) du
= \uint e^{i \lambda t} e_b^\top  \OBM(h)  (i\lambda I_{kq}-\BLAM )^{-1} \BT e_a \Phi_a(d\lambda),
\end{align*}
where $\Phi_a(\cdot)$ is the random spectral measure from  the spectral representation of $\CY_a$. For the second summand, we substitute $Y_a(t)$ as well as its derivatives (cf.~Proposition 2.4 by \citealp{VF23pre}) by their spectral representation.
We obtain 
\begin{align*}
0 &=
\int_{-\infty}^t e_b^\top  \OBM(h)  e^{\BLAM (t-u)} \BT e_a Y_a(u) du 
+ \sum_{m=0}^{p-q-1} e_b^\top  \OBMnu (h) \BT e_a D^{(m)} Y_a(t) \\
&= \uint e^{i \lambda t} e_b^\top  \OBM(h) (i\lambda I_{kq}- \BLAM )^{-1} \BT e_a \Phi_a(d\lambda)\\
&\quad\phantom{h} + \sum_{m=0}^{p-q-1} e_b^\top  \OBMnu (h) \BT e_a \uint (i \lambda)^m e^{i \lambda t} \Phi_a(d\lambda).
\end{align*}
Denoting
$\psi(\lambda,h) = e_b^\top  \OBM(h) (i\lambda I_{kq}-\BLAM )^{-1} \BT e_a
 + \sum_{m=0}^{p-q-1} e_b^\top  \OBMnu (h) \BT e_a (i \lambda)^m,$
for $\lambda \in \R$ and $0\leq h \leq 1$, it follows
that
\begin{align*}
0 = \BE \left[\left| \uint  e^{i \lambda t} \psi(\lambda,h) \Phi_a(d\lambda) \right|^2 \right]
  = \uint \vert \psi(\lambda,h) \vert^2 f_{Y_a Y_a}(\lambda) d\lambda,
\end{align*}
and therefore $\vert \psi(\lambda,h) \vert^2 f_{Y_a Y_a}(\lambda)=0$ for (almost) all $\lambda \in \R$. But by \Cref{Graph well defined}, $f_{Y_aY_a}(\lambda)>0$ for all $\lambda\in \R$, which yields  $ \psi(\lambda,h) =0$
for $0\leq h \leq 1$ and (almost) all $\lambda \in \R$. 
\cite{BE09}, (4.23), provides due to $i\lambda \in \C \setminus \sigma(\BA)$ that
\begin{align*}
(i\lambda I_{kq}-\BLAM)^{-1} = \frac{1}{\chi_\BLAM (i\lambda)} \sum_{j=0}^{kq-1} (i\lambda)^j \Delta_j,
\end{align*}
where $\Delta_j \in M_{kq}(\R)$, $\Delta_{kq-1}=I_{kq}$, and 
$\chi_\BLAM(z) = z^{kq} + \gamma_{kq-1}z^{kq-1}+ \cdots + \gamma_1 z + \gamma_0,$ 
$z\in \C$, is the characteristic polynomial of $\BLAM$ with $\gamma_{1},\ldots,\gamma_{kq-1} \in \R$, $\gamma_{kq}=1$, cf.~\cite{BE09}, (4.4.3). Thus,
\begin{align*}
0 = \psi(\lambda,h) 
= \frac{1}{\chi_\BLAM (i\lambda)} \sum_{j=0}^{kq-1} (i\lambda)^j e_b^\top  \OBM(h) \Delta_j \BT e_a
    +  \sum_{m=0}^{p-q-1} e_b^\top  \OBMnu (h) \BT e_a (i\lambda)^m,
\end{align*}
and multiplication by the characteristic polynomial yields
\begin{align*}
0 
  &=  \sum_{j=0}^{kq-1} (i\lambda)^j e_b^\top  \OBM(h) \Delta_j \BT e_a
    +  \sum_{m=0}^{p-q-1}  \sum_{\ell=0}^{kq} e_b^\top  \OBMnu (h) \BT e_a \gamma_\ell (i\lambda)^{\ell+m}.
\end{align*} 
In the first sum there are powers up to $kq-1$, while in the second sum there are powers up to $kq-1+p-q$. 
For $\ell=kq$ and  $m=0,\ldots,p-q-1$  we receive in the second summand powers higher than $kq$ and their prefactors have to be zero. Due to $\gamma_{kp}=1$
we receive then for $m=0,\ldots,p-q-1$,  
\begin{align*}
e_b^\top  \OBMnu (h) \BT e_a=0.
\end{align*}
Inserting this result into $\psi(\lambda,h)=0$ yields
\begin{align*}
0 
  = e_b^\top  \OBM(h) (i\lambda I_{kq} -\BLAM )^{-1} \BT e_a 
  = \uint e^{-i\lambda s} {\boldsymbol 1}_{ \{s\geq 0 \}} e_b^\top  \OBM(h)  e^{\BLAM s} \BT e_a  ds.
\end{align*}
Together with the already known integrability, 
\cite{PI09}, \mbox{Corollary 2.2.23, provides} 
\begin{align*}
e_b^\top  \OBM(h) e^{\BLAM t} \BT e_a  = 0, \quad t\geq 0,
\end{align*}
which finally concludes the proof of (a).\\
(b)\, Due to the similarity of the results in \Cref{Projektionen für MCARMA} and \Cref{Projektionen mit limneu},
we just have to replace $\BM(h)$ by $\BM'(0)$ and $\BMnu(h)$ by $\BMnu'(0)$ in the proof of (a).
\end{proof}

\begin{proof}[Proof of \Cref{Zweite Charakterisierung gerichtete Kante für MCARMA}] $\mbox{}$\\
(a) \, Based on the characterisations in \Cref{Erste Charakterisierung gerichtete Kante für MCARMA}(a), the same ideas as in the proof of Theorem 6.19(a) in \cite{VF23pre}  can be carried out and therefore, the proof is omitted. First, we replace the matrix exponential $e^{\BA h}$ in \Cref{Erste Charakterisierung gerichtete Kante für MCARMA}(a) by powers of the matrix $\BA$ and second, we replace $e^{\BLAM h}$ by \mbox{powers of $\BLAM$.} \\
(b) \, Follows in analogy to (a) using \Cref{Erste Charakterisierung gerichtete Kante für MCARMA}(b).
\end{proof}

\begin{proof}[Proof of \Cref{Einfache Charakterisierung ungerichtete Kante für MCARMA}] $\mbox{}$\\
(a) \, Based on \Cref{Projektion für MCARMA für S=V}(a), the proof of the first characterisation in (a) can be done in the same way as the proof of  \cite{VF23pre}, Proposition 6.13(a). The second characterisation in (a) follows  along the lines of the proof of \cite{VF23pre}, Theorem 6.19(b).\\
(b) \, Based on \Cref{Darstelung von Yb} and \Cref{Projektion für MCARMA für S=V}(b), statement (b) can be proven  analogously to  \cite{VF23pre}, Proposition 6.13(b).
\end{proof}


\section{Conclusion}
In this paper, we have applied the concept of  (local) orthogonality graphs to state space models. For the state space models, we have assumed that they have a representation in controller canonical form satisfying the mild assumptions of \eqref{Assumption Q} such that there exists a stationary invertible version of the state space model; the term invertible reflects that we are able to recover the state process from the observation process. These assumptions have been summarised under the acronym ICCSS$(p,q)$ model with $p>q>0$. The ICCSS processes satisfy the assumptions of the (local) orthogonality graphs defined in \cite{VF23pre} so that the graphical models are well-defined and several notions of causal Markov properties hold. However, the invertibility of the state process and the representation of the state space model in controller canonical form are not necessary for the existence of (local) orthogonality graphs. The orthogonality graphs exist for a much broader class of state space models, but for the analytic representations of the edges, these additional assumptions are useful. The characterisations of the edges of the ICCSS process require the knowledge of the orthogonal projections of the state process onto linear subspaces generated by subprocesses, and for the derivation of these orthogonal projections the invertibility of the state process is important. The orthogonal projections depend on the model parameters of the controller canonical form and therefore, the edges of the (local) orthogonality graph are also uniquely characterised by these model parameters.

\bibliographystyle{imsart-nameyear}
\bibliography{102_Literatur}
\end{document}